\documentclass{article}
\usepackage{graphicx} 
\usepackage{amsmath,graphicx}

\usepackage{tabularx}  
\usepackage{longtable} 

\usepackage{amsthm}
\usepackage{amssymb}
\usepackage{xypic}
\usepackage{fancyhdr}
\usepackage{multirow}
\usepackage{hyperref}
\usepackage[square,numbers]{natbib}
\theoremstyle{plain}
\newtheorem{theorem}{Theorem}
\newtheorem{proposition}[theorem]{Proposition}
\newtheorem{corollary}[theorem]{Corollary}
\newtheorem{lemma}[theorem]{Lemma}

\theoremstyle{definition}
\newtheorem{definition}{Definition}

\newtheorem{example}{Example}

\newtheorem{question}{Question}
\def\shf{\mathcal}
\def\cshf{\mathfrak}
\def\col{\mathcal}
\def\cat{\bf}

\def\AR{\textrm{AR }}

\def\id{\textrm{id }}
\def\pr{\textrm{pr}}
\def\argmin{\textrm{argmin }}

\newcommand\blfootnote[1]{%
  \begingroup
  \renewcommand\thefootnote{}\footnote{#1}%
  \addtocounter{footnote}{-1}%
  \endgroup
}

\title{Analyzing the topological structure of composite dynamical systems}
\author{Michael Robinson \and Michael L. Szulczewski\blfootnote{Approved for Public Release by The MITRE Corporation; Distribution Unlimited. Public Release Case Number 25-2751.  The author's affiliation with The MITRE Corporation is provided for identification purposes only, and is not intended to convey or imply MITRE's concurrence with, or support for, the positions, opinions, or viewpoints expressed by the author.  ©2025 The MITRE Corporation. ALL RIGHTS RESERVED.} \and James T. Thorson}
\date{September 2025}

\begin{document}

\maketitle

\begin{abstract}
This chapter explores dynamical structural equation models (DSEMs) and their nonlinear generalizations into sheaves of dynamical systems.  
It demonstrates these two disciplines on part of the food web in the Bering Sea.  
The translation from DSEMs to sheaves passes through a formal construction borrowed from electronics called a netlist that specifies how data route through a system.
A sheaf can be considered a formal hypothesis about how variables interact, that then specifies how observations can be tested for consistency, how missing data can be inferred, and how uncertainty about the observations can be quantified.    
Sheaf modeling provides a coherent mathematical framework for studying the interaction of various dynamical subsystems that together determine a larger system.
\end{abstract}

\tableofcontents

\section{Introduction}

Ecologists often study systems on spatial and temporal scales that cannot be experimentally manipulated (ecosystem processes are distributed across continents, and arise from evolutionary dynamics over millennia), and for which extrapolating the results of experiments at fine space-time scales is challenging \cite{tilman_spatial_1997}.  These systems are also challenging to study because observational data can be noisy and sporadic.  A third challenge is the presence of complex, causal relationships between system variables that can change over time.

Understanding the dynamics of these kind of large composite models is much easier reductively.
Roughly speaking, a \emph{subsystem} is a collection of state variables that makes sense as an independent dynamical system (Definition \ref{def:subsystem}).
Subsystems can be isolated for a variety of reasons, in addition to spatial or temporal separation.
Regardless of the reason for the isolation, there is a canonical way to write a dynamical system in terms of its subsystems.
This subsystem decomposition is a convenient way to explore dynamical summaries of the original model (Section \ref{sec:subsystems}).

This chapter explores dynamical structural equation models (DSEMs) and their nonlinear generalizations via a topologically motivated translation into \emph{sheaves of dynamical systems} (Sections \ref{sec:netlist_sheaves} and \ref{sec:subsystems}).  Sheaves are a strict generalization of DSEMs into nonlinear models, which they losslessly represent (Theorem \ref{thm:dsem_solutions_in_sheaf}).
The translation of DSEMs into sheaves follows a clear graphical recipe, which allows handling observations in three  ways: (1) as individual observations, (2) as individual timeseries, and (3) as collections of dynamically related timeseries.

The translation from DSEMs to sheaves passes through a formal construction borrowed from electronics called a \emph{netlist} that specifies how data route through a system.
Because the netlist and sheaf methodology is explicit and graphical, we include several illustrative examples (Figures \ref{fig:linear_regression_netlist} and \ref{fig:linear_regression_sheaf}).  One real-world example involves part of the food web in the Bering Sea (Figure \ref{fig:schematic_sheafify}; Sections \ref{sec:bering_dsem}, \ref{sec:bering_sheaf}, and \ref{sec:bering_sheaf_subsystems}).

Sheaves provide many advantages to a modeler.  They enable exploring the impact of uncertainty in various ways.  They support inference of missing or erroneous data, including system parameters and coefficients (Section \ref{sec:netlist_sheaves}).  They also enable forecasts and retrocasts through the same ``interface,'' namely \emph{consistency radius optimization} (Section \ref{sec:bering_sheaf}). 

Sheaves also highlight the importance of the original DSEM in model summarization.
Using the sheaf of subsystems, Corollary \ref{cor:acyclic_dsem_subsystems} shows that the subsystems of a DSEM can be ``read off'' its associated graph.  This is applied to the Bering Sea ecosystem model in Section \ref{sec:bering_sheaf_subsystems}.

\subsection{Related work}

The challenges in modeling ecological systems have motivated interest in structural causal models (SCMs) \cite{pearl_causal_2009}.  SCMs can be fit to observational data in space and time, and can decompose the total effect of one variable on another via a combination of direct and indirect effects \cite{grace_scientists_2020,arif_applying_2023}.  
Recently, SCMs have been adapted to the analysis of ecological time series via DSEMs \cite{thorson_dynamic_2024}. 

The key idea behind SCMs is that systems can be understood by decomposing them into coherent subsystems.  The idea of reducing systems into subsystems has a long history, with general mathematical descriptions of composite systems given by the field of cybernetics, for which \citet{joslyn_2003} and \citet{Ashby_1956} are good introductions.
Beyond cybernetics, the study of subsystems of dynamical models \cite{Yin2016SubsystemDF} has occurred in many fields, including manufacturing and operations research \cite{Wynn2018ProcessMI,Suh2015SeeingCS,Johnson1993SubsystemDI}, design \cite{Alexander_1964}, statistical physics \cite{Zwanzig_2001}, mathematical systems \cite{Chorin_2002}, biology \cite{Liberzon2015TheMS}, and chemistry \cite{Hirono_2021}.  

Although algorithmic and systematic decomposition of systems into subsystems have become common since the dawn of cybernetics, it remains challenging.  \citet{Maier2017ModelGI} laments, ``Even though abstraction is frequently mentioned with regards to modeling and simulation, formal definitions are harder to find.''  One challenge is that decompositions are often not unique: for example, one may choose to group state variables based on constraints rather than functional units \cite{Chiriac2011ThreeAT,Komoto_2011}. 
These choices are important because they drive the usefulness of the decomposition~\cite{Maier2017ModelGI}.  
For example, overlapping, rather than disjoint, subsystem decompositions are useful for analyzing stability of an entire system \cite{Sloth_2016,Anderson_2010}.

We argue that a properly general and formal definition of a subsystem decomposition must support overlappingness, non-uniqueness, and ambiguous granularity.  
Because the collection of all subsystems forms a \emph{mathematical sheaf} (Definition \ref{df:subsystem_sheaf}), this implies that seeking disjoint, unambiguous subsystems (as is often done) is fraught.  

Aspects of the formalism we introduce in this chapter are not entirely novel.  For instance, \citet{Hirono_2021} defines a \emph{CRN morphism} that is a special case of our Definition \ref{def:subsystem}.  Additionally, the sheaf of subsystems is based upon a clear graphical representation, which is well known in the analysis of software \cite{Montes_1998,Abadi_1993}.  Moreover, \citet{Abadi_1993} uses the term \emph{refinement mapping}, which evokes the analogous term from sheaves (Definition \ref{df:presheaf}).

Roughly dual to the notion of a subsystem is that of an \emph{invariant set} of a dynamical system (our Definition \ref{def:subsystem} makes this a \emph{true} duality).  Invariant sets are widely used in dynamical systems \cite{Strogatz}, where they generalize equilibrium sets and attractors.  For linear systems, duality between invariant sets and subsystems is immediate and useful.  For instance, the design structure matrix \cite{Steward_1981} yields invariant sets, giving a clear duality to subsystems.

Finally, we note that the discipline of modeling a system's state via a decomposition into subsystems of state equations is explained in detail in \citet[Sec. 5]{Robinson_multimodel}, and is specialized to subsystem graphs in \citet{Kearney_2020}.  In \citet{Kearney_2020}, the dynamics are specified locally and are much easier to specify due to the fact that the system is given a graph structure.

\subsection{Contributions}

This chapter provides an introduction to the discipline of modeling and analyzing a composite system using the language and tools of topology, centered around \emph{sheaves}.
Sheaf modeling provides a coherent mathematical framework for studying the complicated interaction of various dynamical subsystems that together determine a larger system.
The guiding principles of sheaf modeling are that
\begin{itemize}
    \item a \emph{sheaf} represents a hypothesis about how variables will interact (Definition \ref{df:sheaf}),
    \item a non-global \emph{assignment} represents the observations collected on the variables in its support (Definition \ref{def:assignment}),
    \item \emph{minimizing consistency radius} estimates values of the variables and parameters that were not observed (Definition \ref{df:consistency_radius}), and
    \item the \emph{minimal consistency radius} is a measure of the consistency between the observations and the hypothesis.
\end{itemize}
This chapter shows that when a dynamical system is described by a linear system, there are three sheaves that provide increasingly granular data about the interactions between variables:
\begin{enumerate}
    \item the \emph{sheaf of subsystems} (Definition \ref{df:subsystem_sheaf}),
    \item the \emph{netlist sheaf} with timeseries as stalks (Definition \ref{df:netlist_sheaf}), and
    \item the \emph{netlist sheaf} with additional stalks for individual observations (Definition \ref{df:netlist_sheaf_individual}).    
\end{enumerate}

\subsection{Chapter outline}

Section \ref{sec:ecology} describes a model of a food web in the Bering Sea, which we use to illustrate the use of sheaves.  This system is large enough to exhibit interesting structures, and corresponding observational data \cite{thorson_dynamic_2024} are available.  Additionally, we present a graphical causal modeling discipline called \emph{dynamical structural equation modeling} that serves as an entry point into the more sophisticated (but admittedly less familiar) topological sheaf models.  As is later shown in Section \ref{sec:netlist_sheaves}, sheaves are a strict generalization of DSEMs.  Sheaves can be nonlinear, whereas DSEMs are linear.

Section \ref{sec:netlist_sheaves} constructs sheaves that model composite systems, and develops the main inferential tool, \emph{consistency radius minimization}.  Section \ref{sec:netlist_sheaves} is self-contained, as all of the mathematical background necessary to understand the constructions is introduced as it is needed.  Small concrete examples of the construction and use of sheaf models are presented to build intuition as well.  

In Section \ref{sec:bering_sheaf}, we revisit the ecological model from Section \ref{sec:ecology} using the sheaf tools from Section \ref{sec:netlist_sheaves}.  The interface between observational data, sheaves, and their inference tools is explored in detail.  Moreover, we compare differences between the DSEM and sheaf approaches in detail.

Section \ref{sec:subsystems} introduces the idea of a general topological \emph{dynamical system}, and shows that every dynamical system induces a \emph{sheaf of subsystems} and a \emph{cosheaf of invariant sets}, which form a dual pair.  We prove that under appropriate conditions, the subsystems of a DSEM can be ``read off'' rather directly (Corollary \ref{cor:acyclic_dsem_subsystems}).  This provides theoretical justification for why DSEMs are a useful way to describe a composite linear system by way of its subsystems. 

Section \ref{sec:bering_sheaf_subsystems} revisits the ecological model from Section \ref{sec:ecology} once again.  Because the model satisfies the hypothesis of Corollary \ref{cor:acyclic_dsem_subsystems}, we are able to present a clear representation of all the subsystems present in the model.

Finally, Section \ref{sec:conclusion} concludes the chapter with practical advice for modelers and a brief discussion of future research work. 

\section{Dynamical modeling of ecosystems}
\label{sec:ecology}

This section begins with a brief recount of modeling linear dynamical systems according to an underlying graph structure, and then presents a representative ecosystem model that will be revisited several times in the chapter.

\subsection{DSEM background and motivation}
\label{sec:dsem_background}

\begin{definition}
\label{def:dsem}
    Given a set of \emph{variables} $X =\{x_1, \dotsc, x_J\}$, and a set $Y =\{ t_1 < \dotsb < t_T\}$ of real valued \emph{time lags},
    a \emph{dynamic structural equation model (DSEM)} consists of an edge-labeled directed graph $G$ with vertices $X \times Y$ and edges $E$ such that 
    \begin{description}
        \item[Causality] The presence of an edge $(x_{j_1},t_{k_1}) \to (x_{j_2},t_{k_2})$ implies that $t_{k_1} \le t_{k_2}$, and
        \item[Linearity] Each edge $(x_{j_1},t_{k_1}) \to (x_{j_2},t_{k_2})$ is labeled with a real number $\gamma_{j_1,k_1,j_2,k_2}$ called the \emph{path coefficient} for that edge.
    \end{description}
    The absence of an edge in the graph is assumed to be equivalent to assigning a path coefficient of $0$.
    For brevity, we write a vertex $(x_j,t_k)$ simply as $x_{j,k}$.
\end{definition}

The variables in a DSEM are to be interpreted as $C^1(\mathbb{R})$ functions, which are continuous timeseries.
A directed edge $x_{i,j} \to x_{i',j'}$ is to be interpreted as specifying that a change in $x_{i}$ causes a proportional (linear) change in $x_{i'}$ after a lag of $(t_{j'} - t_{j})$, with magnitude controlled by the associated path coefficient $\gamma_{i,j,i',j'}$.
Under this interpretation, a DSEM implies that a first order system of linear differential equations governs the values of the variables:
\begin{equation}
\label{eq:dsem_exact}
    \frac{dx_k(\tau - t_\ell)}{d\tau} = \sum_{i=1}^J\sum_{j=1}^T \gamma_{k,\ell,i,j} x_{i}(\tau - t_j).
\end{equation}
In what follows, we will refer to solutions of Equation \ref{eq:dsem_exact} as \emph{solutions to the DSEM}.

In the use of Equation \eqref{eq:dsem_exact} with observational data, there are two kinds of errors that need to be considered: \emph{exogenous} errors and \emph{measurement} errors.  Exogenous errors accumulate, which means that an error in the value of a variable $x_k$ at given time $\tau$ impacts the value of $x_k$ at all later times.  As a result, there is a dependence between the exogenous errors of $x_k$ at different times.  In contrast, measurement errors at different times are assumed to be independent.

Exogenous errors will be represented by an additive term, $\epsilon_{k,\ell}$, resulting in
\begin{equation}
\label{eq:dsem_errored}
    \frac{dx_k(\tau - t_\ell)}{d\tau} = \sum_{i=1}^J\sum_{j=1}^T \gamma_{k,\ell,i,j} x_{i}(\tau - t_j) + \epsilon_{k,\ell}(\tau).
\end{equation}

We can approximate the solution to Equation \eqref{eq:dsem_errored} using the one-step backwards Euler method with time step $h$,
\begin{equation*}
    \frac{dx_k(\tau - t_\ell)}{d\tau} \approx\frac{1}{h}\left(x_k(\tau-t_{\ell}) - x_k(\tau-t_{\ell} - h)\right),
\end{equation*}
so that Equation \eqref{eq:dsem_errored} becomes a system of $M=TJ$ linear algebraic equations,
\begin{equation}
\label{eq:dsem_errored_approx}
    x_k(\tau-t_{\ell}) \approx x_k(\tau-t_{\ell} - h) + h \sum_{i=1}^J\sum_{j=1}^T \gamma_{k,\ell,i,j} x_{i}(\tau - t_j) + h \epsilon_{k,\ell}(\tau).\\
\end{equation}

If we fix a value of $\tau$ and organize the set of values $\{x_k(\tau-t_{\ell})\}$ into a vector $\mathbf{X}$ of length $M$), Equation \eqref{eq:dsem_errored_approx} can be compactly written in matrix form as
\begin{equation}
\label{eq:dsem_matrix}
    \mathbf{X} \approx \mathbf{P} \mathbf{X} + \mathbf{E},
\end{equation}
where the entries of the $M \times M$ \emph{path coefficient matrix} $\mathbf{P}$ contain both the path coefficients from the DSEM (scaled by $h$) and the additional nonzero entries due the $x_k(\tau-t_{\ell} - h)$ terms.
In what follows, we will take $h=1$, so that the path coefficients in the DSEM appear unchanged as elements of the matrix $\mathbf{P}$.

To obtain the path coefficient matrix $\mathbf{P}$ from observations of $\mathbf{X}$, we assume the exogenous errors follow a multivariate normal distribution with variance $\mathbf{V}$, namely
\begin{equation*}
    \mathbf{E} \sim \text{MVN}(\mathbf{0,V}),
\end{equation*}
where $\mathbf{E}$ is the length $M$ vector containing errors $\epsilon_{tj}$.  

Equation \eqref{eq:dsem_matrix} can then be re-arranged to yield a Gaussian Markov random field,
\begin{align}
    \mathbf{X} &\sim \text{MVN}(\mathbf{0,Q}^{-1}) \\
    \mathbf{Q} &= (\mathbf{\id-P}^T) \mathbf{V}^{-1} (\mathbf{\id-P}),
\end{align}
where $\id$ is the identity matrix.  The path coefficient matrix $\mathbf{P}$ can be obtained from the Cholesky decomposition of $\mathbf{Q}$.
The necessary calculations can be efficiently evaluated using sparse libraries, such as {\tt Eigen} and {\tt CHOLMOD} \cite{davis_university_2011}, and we use Template Model Builder \cite{kristensen_tmb:_2016} to incorporate automatic differentiation and implement the Laplace approximation \cite{skaug_automatic_2006} to marginalize across random effects.

Now we address measurement errors.  
Assume the distribution of measurement errors of the variable $x_k$ is given by a distribution $f_j$ parameterized by $\theta_j$ at time $t_j$.  (If one does not wish to model measurement errors explicitly, so that measurement errors are entirely captured by the exogenous error term, this is obtained by choosing $f_j$ so that it has probability $1$ at $x_{k,j}$.) 
Let us write $y_{k,j}$ for the observation of the variable $x_{k,j}$.
We therefore can express the mean of the distribution of $y_{k,j}$ through a link function $g_j$, via
\begin{equation*}
    y_{k,j} \sim f_j \left( g_j^{-1}( \mu_j + x_{k,j}), \theta_j \right),  
\end{equation*}
where $\mu_j$ is the true mean.

The clearest way to obtain the required sparsity in solving for $\mathbf{P}$ is to assume additionally that the measurement errors for a given variable do not depend on time $t_j$.
Let $\mathbf{G}$ be the $J\times J$ matrix that is diagonal, and whose diagonal terms are given by the link functions $g_j$.
With this in hand, $\mathbf{V}$ takes the form
\begin{align}
  \mathbf{V} = \id_{T\times T} \otimes \mathbf{GG}^T,
\end{align}
where $\otimes$ is the Kronecker product.
This implies that $\mathbf{V}$ is block diagonal, and is thereby efficient to invert.

\subsection{Ecological background and the DSEM system for the Bering Sea}
\label{sec:bering_dsem}

To demonstrate the use of sheaves for dynamical systems, we make a sheaf from a DSEM for ecological mechanisms linking regional oceanography (winter sea ice extent) to first-winter survival of juvenile Alaska pollock (\textit{Gadus chalcogrammus}) in the eastern and northern Bering Sea \cite{thorson_dynamic_2024}.  The model starts by specifying that abundance of age-0 pollock $R_t$ (termed ``age-0 recruitment") can be predicted from the biomass of spawning females $S_t$ in a given year $t$:

\begin{equation}
    R_t = S_t e^{\alpha-\beta S_t + \epsilon_t}
\end{equation}
where $e^\alpha$ is the maximum expected recruits per spawning biomass, $\beta$ is the expected density-dependent decrease in recruits per spawning biomass as biomass increases, and $\epsilon_t$ is additional process error representing unmodeled variation in recruitment.  This ``Ricker stock-recruit model" \cite{ricker_stock_1954} has been used for over 70 years to represent density-dependent changes in juvenile survival, and as the basis for defining biological reference points that are used worldwide to identify sustainable levels of fishing mortality \cite{smith_gospel_2001}.  The Ricker model is expected to arise for species where adult abundance directly impacts juvenile survival---for example, due to cannibalism or interference competition \cite{foss-grant_hierarchical_2016}.  Alaska pollock are cannibalistic, so the Ricker model has theoretical justification.  Usefully, the Ricker model can be linearized as:
\begin{equation}
    \text{log} \left( \frac{R_t}{S_t}\right) = \alpha - \beta S_t + \epsilon_t
\end{equation}
and a DSEM can be used to elaborate the mechanisms that contribute to process errors $\epsilon_t$ based on prior ecological hypotheses.  

The DSEM we translate into a sheaf was previously developed by~\citet{thorson_dynamic_2024}.  It specifies that variable winter sea ice formation ($\mathit{SeaIce}$) drives residual variation in log-recruits per spawning biomass ($\mathit{Survival}$) via two paths, mediated by sea-ice impacts on either copepod abundance ($\mathit{Copepod}$) or krill abundance ($\mathit{Krill}$), and resulting consumption by juvenile pollock.  See Table \ref{tab:bering_sea_variables} and \ref{tab:bering_sea_hypotheses} for more details on the variables and mechanisms in the model.  The DSEM includes a first-order autoregressive term for each variable, to allow the model to correct for bias that can arise when correlating variables that follow an autoregressive process (summarized in \cite{mccallum_is_2010}).  This first-order autoregression can also be interpreted to represent Gompertz density-dependence and therefore has some scientific interest \cite{knape_are_2012}, although it is not further discussed here. 

\begin{table}
  \caption{Variables that describe Alaska pollock recruitment used in the DSEM and sheaf.  All except $\mathit{Spawners}$ are transformed by the natural logarithm and then centered (i.e., subtracted by their mean) prior to analysis.  Timeseries of the variables are taken from \cite{thorson_dynamic_2024}.}
\begin{center}
\begin{tabularx}{\textwidth}{| l |  X | } 
  \hline
  Name & Description \\ 
  \hline

  $\mathit{SeaIce}$ & Average spatial extent (km$^2$) of sea ice in the Bering Sea from Oct.15 to Dec.15 the preceding year, from the National Snow and Ice Center's Sea Ice Index, Version 3 \cite{fetterer_sea_2017} \\ & \\

  $\mathit{ColdPool}$ & Spatial extent (km$^2$) of waters with temperatures $\leq 2^\circ$C near the seafloor, interpolated from measurements by the eastern Bering Sea bottom trawl survey and compiled in R-package ``coldpool" \cite{rohan_evaluating_2022} \\ &\\

  $\mathit{Spawners}$ & Female spawning biomass (in units of $10^6$ kg) for Alaska pollock in the eastern and northern Bering Sea, estimated by the age-structured stock assessment model used for management \cite{ianelli_assessment_2022} \\ &\\

  $\mathit{Survival}$ & Age-0 recruits per spawning biomass ($10^3$ count/kg), calculated as age-1 abundance the following year  ($10^9$ count) estimated by the age-structured stock assessment model \cite{ianelli_assessment_2022} divided by $\mathit{Spawners}$ \\ &\\

  $\mathit{Copepods}$ & Density of $\geq 2$ mm copepods (count/m$^3$) from the Bering Sea middle shelf \cite{siddon_ecosystem_2022}, averaged across samples obtained during the fall mooring cruise along the $70 \; $ isobath from Sept. to early Oct. \cite{duffy-anderson_return_2017} (calculated by Dave Kimmel, pers. comm.) \\ &\\

  $\mathit{Krill}$ & Index of euphausiid abundance (count/m$^3$) \cite{ressler_developing_2012} obtained from backscatter measured during a summer acoustic-trawl survey in the eastern Bering Sea and converted to abundance using a target-strength model \cite{smith_distorted_2013} \\& \\ 

  $\mathit{DietCopepods}$ & Biomass of copepods divided by total prey biomass in juvenile stomach samples (kg/kg), calculated from a fall surface-trawl survey in the eastern Bering Sea \cite{murphy_northern_2021}.  For each surface trawl, total catch of juvenile pollock is weighed, individual pollock are subsampled, and stomach contents for subsampled individuals are identified to species and weighed. The diet index is calculated as the average across subsampled stomachs, weighted by the catch of juvenile pollock in the associated surface trawl sample (calculated by Alex Andrews, pers. comm.).  \\ &\\

  $\mathit{DietKrill}$ & Same as $\mathit{DietCopepods}$, but for euphausiids (krill) \\

  \hline
\end{tabularx}
  \label{tab:bering_sea_variables}
\end{center}
\end{table}

\begin{table}
  \caption[Hypotheses for Alaska pollock recruitment]{List of path coefficients connecting variables (defined in Table \ref{tab:bering_sea_variables}), supporting ecological hypotheses, and hypothesized sign for the path used in the DSEM case study. We also include a first-order autoregressive term for each variable (i.e., $8$ AR1 coefficients, not shown here) for reasons discussed in Section \ref{sec:bering_dsem}.}
\begin{center}
\begin{tabularx}{\textwidth}{| l | X | c | } 
  \hline
  Path & Ecological hypothesis and evidence & Sign \\ 
  \hline

  $SeaIce \rightarrow ColdPool$ & Sea ice formation ($\mathit{SeaIce}$) causes variation in summer cold-pool extent ($\mathit{ColdPool}$) & $+$ \\ & & \\ 

  $ColdPool \rightarrow Copepods$ & Warmer water temperatures ($\mathit{ColdPool}$) result in higher copepod metabolism and therefore earlier onset of winter diapause, resulting in a decrease in fall copepod abundance ($\mathit{Copepods}$) \cite{coyle_calanus_2017} & $+$ \\ & & \\

  $ColdPool \rightarrow Krill$ & Water temperatures ($\mathit{ColdPool}$) might affect krill overwinter survival, affecting summer krill abundance ($\mathit{Krill}$) & $?$ \\ & & \\  

  $Copepods \rightarrow DietCopepods$ & Increased copepod abundance will result in them being a higher proportion of age-0 fall stomach contents ($\mathit{DietCopepods}$), due to pollock being hypothesized to be a relative non-selective predator & $+$ \\ & & \\  

  $Krill \rightarrow DietKrill$ & Same as $Copepods \rightarrow DietCopepods$ but for krill & $+$ \\ & & \\  

  $DietCopepods \rightarrow Survival$ & Increased fraction of fall diet from copepods ($\mathit{Copepods}$) will increase energy reserves and subsequent survival of age-0 over their first winter ($\mathit{Survival}$) \cite{hunt_climate_2011} & $+$ \\ & & \\  

  $DietKrill \rightarrow Survival$ & Same as $DietCopepods \rightarrow Survival$, but for krill  & $+$ \\ & & \\  

  $Spawners \rightarrow Survival$ & Increased spawning biomass ($\mathit{Spawners}$) will cause a density-dependent decrease in survival ($\mathit{Survival}$) \cite{foss-grant_hierarchical_2016} & $-$ \\ & & \\  

  \hline
\end{tabularx}
  \label{tab:bering_sea_hypotheses}
\end{center}
\end{table}

\section{Sheaf encodings of composite systems}
\label{sec:netlist_sheaves}

In this section, we explain how to construct a \emph{netlist sheaf} whose \emph{global sections} correspond bijectively to the solutions of a DSEM.
This is performed in two main steps: (1) the DSEM is translated into a \emph{netlist}, and (2) the netlist is translated into the \emph{netlist sheaf}.
Since the machinery of sheaves is not in wide usage, Section \ref{sec:sheaf_background} provides the necessary background.

With the machinery and the translation in place, Theorem \ref{thm:dsem_solutions_in_sheaf} establishes that the two representations, the DSEM and the netlist sheaf, are equivalent.  
The \emph{global sections} of the netlist sheaf are in bijective correspondence with solutions to the DSEM.
Moreover, a process called \emph{consistency radius minimization} in the sheaf finds approximate solutions to the DSEM, and this process is robust to perturbations.

\begin{figure}
    \centering
    \includegraphics[width=\linewidth]{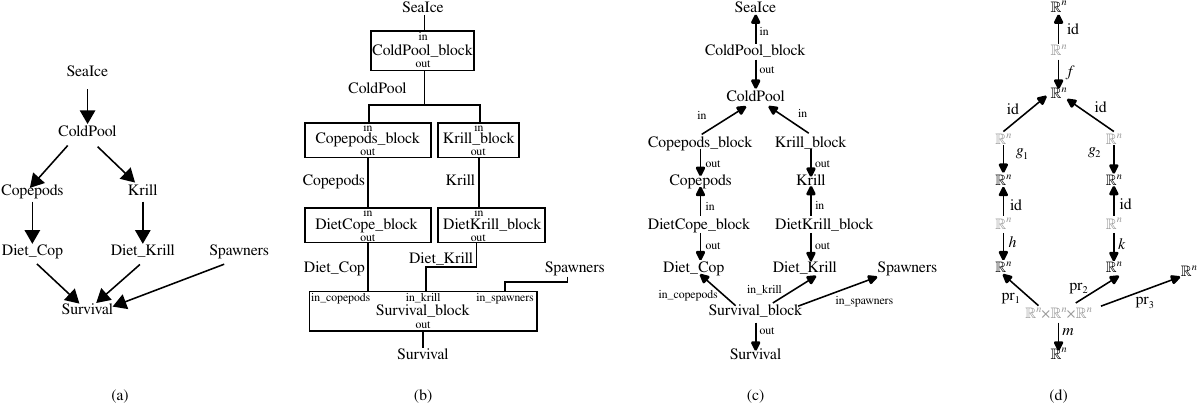}
    \caption{(a) The DSEM model for part of a food web in the Bering Sea \cite{thorson2024dynamic}, (b) its wiring hypergraph, (c) its netlist graph, and (d) its sheaf diagram.  The arrows in each subfigure have different meanings:  in (a) they denote causal, linear relationships (Sec.~\ref{sec:dsem_background}); in (c), they point from netlist parts to nets (Sec.~\ref{sec:netlists}); and in (d), they denote restriction functions (Sec.~\ref{sec:sheaf_background}).  While the DSEM also estimates a first-order autoregressive term for each variable (not shown in (a) to simplify presentation), there is no autoregressive structure assumed in the sheaf model.  This remedied in Section \ref{sec:sheaf_autoregression}.}
    \label{fig:schematic_sheafify}
\end{figure}

Throughout this section, we refer to Figure \ref{fig:schematic_sheafify} for intuition. Figure \ref{fig:schematic_sheafify}(a) shows the DSEM for part of the food web in the Bering Sea.  The DSEM-to-netlist translation, described in Section \ref{sec:netlists}, results in Figure \ref{fig:schematic_sheafify}(b).
Figure \ref{fig:schematic_sheafify}(c) shows a different representation of the netlist that is more expedient for the construction of the netlist sheaf.
Proposition \ref{prop:netlist_hypergraph_and_graph} establishes that the two representations of netlists (Figures \ref{fig:schematic_sheafify}(b)--(c)) determine each other, 
so we may use whichever is more convenient.
Finally, the netlist-to-sheaf translation, described in Section \ref{sec:netlist_sheaves}, results in Figure \ref{fig:schematic_sheafify}(d).
Section \ref{sec:sheaf_autoregression} shows how to encode autoregressive timeseries models as netlist sheaves, which ultimately makes handling missing data both transparent and automatic within the netlist sheaf. 

\subsection{Netlists}
\label{sec:netlists}

The term ``netlist'' appears to have entered the technical lexicon in the early days of computing, when IBM started to automate the wiring of mainframe back planes \cite{ibm_netlists}.  Since that time, the term ``netlist'' has been in wide usage but often without a precise definition.  In order to formalize the concept, we say that a \emph{netlist} describes a system of \emph{parts} interconnected with \emph{nets}, which carry \emph{time-varying signals} (briefly, \emph{variables}).

Each \emph{variable} consists of the specification of a set of possible \emph{values} for a net.
In this chapter, the values for a variable in a net are initially assumed to be continuous timeseries, usually of the form $C^1(\mathbb{R})$.
We will also consider sampled timeseries of the form $\mathbb{R}^{n}$, where $n$ is the length of the timeseries.
In Section \ref{sec:sheaf_autoregression}, we show how to handle missing values in such a timeseries.

Each part has a number of \emph{ports}, to which connections can be made.
Each \emph{port} is either
an \emph{output}, which means that it determines the value of the variable of a net connected to it, or
an \emph{input}, which means that it does not determine the value of the variable of a net connected to it.

Each \emph{net} specifies that a collection of distinct ports on a pair of parts (which need not be distinct) are connected, 
with the requirement that not more than one of these ports be an output.
Finally, each part specifies an \emph{input-output function} for each output port.
The domain of an input-output function is from the product of the set of its input variables,
and its codomain (range) is the set of output variables at the output port.

This formulation leaves open the possibility of nets that are not attached to any output ports, which are called \emph{external inputs}, 
and nets which are not attached to any input ports, which are called \emph{external outputs}.
Clearly each external output must attach to exactly one port, which must be an output port.

\begin{figure}
    \centering
    \includegraphics[width=0.5\linewidth]{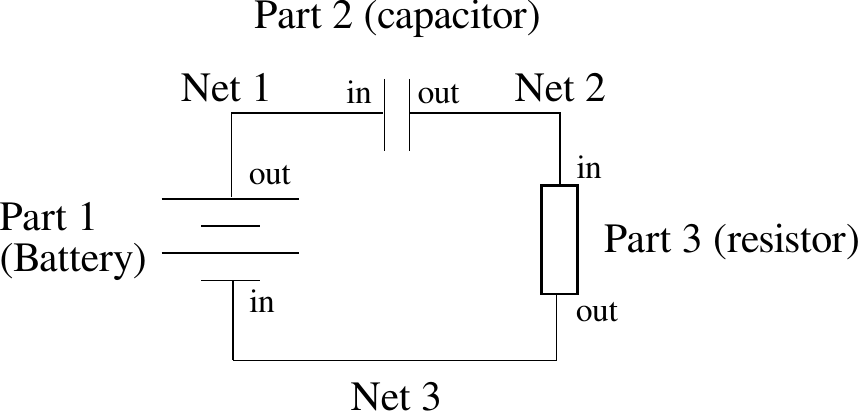}
    \caption{A netlist for an electric circuit, described in Example \ref{eg:eg_netlist}.}
    \label{fig:eg_netlist}
\end{figure}

\begin{example}
    \label{eg:eg_netlist}
    Figure \ref{fig:eg_netlist} shows an electrical circuit with three parts: a battery, a capacitor, and a resistor.  
    These parts are connected to each other by three nets: 
    \begin{enumerate}
        \item Connecting the positive (output) port of the battery to the input port of the capacitor, 
        \item Connecting the output port of the capacitor to the input port of the resistor, and
        \item Connecting the output port of the resistor to the input port of the battery.
    \end{enumerate}
    The values of the variables on the nets specify electrical currents flowing along them.
    We note that the labeling ports as ``input'' and ``output'' in this kind of circuit is arbitrary, since the electrical current can flow in either direction along a net.
    The input-output functions simply recount classical Ohm's law for each of the parts in the circuit.
    This circuit contains no external inputs nor external outputs. 
\end{example}

A DSEM graph can be translated into a netlist via the following construction.

\begin{definition}
\label{df:netlist_from_dsem}
    Given a DSEM, its corresponding netlist is given by the following recipe:
    \begin{itemize}
    \item each DSEM variable (node) becomes a net,
 \item each DSEM variable with more than one input becomes a part,
 \item each net is connected to input ports via its out-neighbors,
 \item each net is connected to output ports via matching the name of the net to the part with the same name (if any exist), and
 \item the part's input-output function is collected from the matrix block in Equation \eqref{eq:dsem_matrix} corresponding to the input and output variables.
    \end{itemize}
\end{definition}

There are two combinatorial structures associated to a netlist, the \emph{wiring hypergraph} and the \emph{netlist graph}.

\begin{definition}
    The \emph{wiring hypergraph} of a netlist is a vertex- and edge-labeled partition-directed multi-hypergraph that has a vertex for each part and an hyperedge for each net.
    
    The label on each vertex is simply the name of the part corresponding to that vertex.
    
    The vertices within a hyperedge correspond to the parts connected to the corresponding net.
    The label on each hyperedge is an ordered triple, consisting of the inputs port of the net (if any), the output port of the net (if any), and the variable name of the net.
    The partition direction of each hyperedge separates the output port from the input ports; either of these may be empty.
\end{definition}

Because the labeling on the wiring hypergraph is complicated, we represent it with a standard visual grammar borrowed from electronics.
Each part is represented by a rectangle with its label in the center of the rectangle.
Each net is drawn as a path (with right-angle bends as needed) to connect the corresponding parts.
If a net has more than two ports, the path is drawn as a tree structure.
The label of the variable of the net is shown next to the path,
but the name of the net's input and output ports are shown inside the connected parts' rectangles, around the edge of the rectangle.
The input-output functions are not shown explicitly.

Figure \ref{fig:schematic_sheafify}(b) shows the wiring hypergraph for the netlist constructed using Definition \ref{df:netlist_from_dsem} for the Bering Sea DSEM.
Notice that the net $\mathit{ColdPool}$ corresponds to a hyperedge of size $3$ in the wiring hypergraph, because it is connected to one output port and two input ports.

\begin{proposition}
\label{prop:dsem_netlist}
    The solutions to a DSEM are in bijective correspondence with labelings of the nets with values of variables that are consistent with the netlist's input-output functions.
\end{proposition}
\begin{proof}
    The solutions to the DSEM are characterized by Equation \eqref{eq:dsem_matrix}, which is a matrix block assembly of everything that is needed to construct the netlist.

    Assume we have a set of variables for all nets that are consistent with the input-output functions.
    As noted above, each variable takes values in a set of the form $C^1(\mathbb{R})$.  
    On the other hand, each input-output function was constructed from a matrix block in Equation \eqref{eq:dsem_matrix}.
    Because all of the DSEM variables appear as nets in the netlist, all such matrix blocks appear as input-output functions somewhere in the netlist.
    This means that Equation \eqref{eq:dsem_matrix} is satisfied by construction.

    Assume that we have a solution to Equation \eqref{eq:dsem_matrix}.
    Definition \ref{df:netlist_from_dsem} constructed the input-output function from the subblock of Equation \eqref{eq:dsem_matrix}, so there is nothing further to prove.
\end{proof}

The wiring hypergraph is closely related to the DSEM, but
for constructing the netlist sheaf in Section \ref{sec:netlist_sheaves}, it is more convenient to use another combinatorial representation.

\begin{definition}
    The \emph{netlist graph} is a vertex- and edge-labeled directed graph that has a vertex for each part, a vertex for each variable, and two edges for each net.
    The label on a vertex is simply the name of the corresponding part or variable.
    The two edges for each net are defined as follows. 
    The first edge is labeled with the input port of the net, and leads from that corresponding part to the net.
    The second edge is labeled with the output port of the net, and leads from that corresponding part to the net.
\end{definition}

Figure \ref{fig:schematic_sheafify}(c) shows the netlist graph for the Bering Sea example.

\begin{corollary}
    The netlist graph is a directed acyclic graph, and induces a preorder on the set of parts and variables.
    In the preorder, each variable is above the parts to which it is connected.
\end{corollary}

\begin{proposition}
\label{prop:netlist_hypergraph_and_graph}
    The netlist graph is the incidence bipartite graph of the wiring hypergraph,
    whose edges are labeled by projecting out the first and second components of the labels of the hyperedges.
Consequently, the netlist graph and the wiring hypergraph determine each other fully.
\end{proposition}

As we will see, the correspondence between the wiring hypergraph and the netlist graph is convenient.
Although Proposition \ref{prop:dsem_netlist} showed that the wiring hypergraph is most closely related to the DSEM, 
we will later show that the netlist graph is most closely related to the netlist sheaf (Theorem \ref{thm:dsem_solutions_in_sheaf}).

\subsection{Sheaves and cosheaves}
\label{sec:sheaf_background}

Sheaves and cosheaves are topological constructions that allow one to study the local consistency structure of a model.
In the case of a DSEM, locality is useful because variables that are near one another in the graph are likely to be related.
This nearness can be most easily formalized by using the netlist graph defined in the previous section.

Since the netlist graph is a directed acyclic graph, it naturally induces a pre-ordered set on the vertices.
That is, if $a \to b$ in a directed graph, we define $a \le b$.  
When the graph is directed and acyclic, generalizing $\le$ to paths within the graph results in a relation $\le$ that is reflexive and transitive.
Pre-ordered sets have a natural notion of neighborhoods, hence a natural topology.

A \emph{topological space} is a mathematical formalism that captures the notion of ``neighborhoods.''  

\begin{definition}
A \emph{topology} on an arbitrary set $X$ is a collection $\col{T}$ of subsets of $X$ satisfying the following four axioms:
\begin{description}
    \item[Empty set] The empty set $\emptyset$ is an element of $\col{T}$,
    \item[Whole set] The set $X$ is an element of $\col{T}$,
    \item[Finite intersection] If $U$ and $V$ are elements of $\col{T}$, then $U \cap V$ is an element of $\col{T}$, and
\item[Arbitrary union] If $\col{U} \subseteq \col{T}$ then $\cup \col{U}$ is an element of $\col{T}$.
\end{description}
The ordered pair $(X,\col{T})$ is called a \emph{topological space}.

Often, rather than specifying $\col{T}$ directly, we specify a collection of subsets $\col{U}$ of $X$ that \emph{generate} the topology, which is the smallest topology (in the sense of inclusion) that contains $\col{U}$.
\end{definition}

The following are elementary examples of topological spaces,
\begin{description}
    \item[Discrete topology] For any set $X$, let $\col{T}$ be the power set of $X$,
    \item[Trivial topology] For any set $X$, let $\col{T} = \{\emptyset, X\}$,
    \item[Euclidean topology] For $X = \mathbb{R}$, the usual topology $\col{T}$ is generated by the set of open intervals $(a,b)$ for $a < b \in \mathbb{R}$.
\end{description}

Additionally, there is a powerful combinatorial theory of topological spaces $(X,\col{T})$ in which the topology $\col{T}$ is a finite set \cite{Barmak_2011}.
For our purposes, the most interesting of these \emph{finite topological spaces} are those that arise naturally from a pre-ordered set, given by the definition below.

\begin{definition}
\label{df:alexandrov}
  Suppose that $(P,\le)$ is a pre-ordered set, which is to say that $\le$ is a reflexive and transitive relation.  The \emph{Alexandrov topology $Alex(P,\le)$ on $(P,\le)$} is the topology generated by all subsets of $P$ of the form $U_x = \{x \le y : y \in P\}$.
\end{definition}

The idea of sheaves and cosheaves is that each \emph{open set}---an element of the a topology---is associated with a set of values, called the stalk (for sheaves) or costalk (for cosheaves). 

\begin{definition}
  \label{df:presheaf}
Suppose $(X,\col{T})$ is a topological space.  A \emph{presheaf $\shf{S}$ of sets on $(X,\col{T})$} consists of the following specification:
\begin{enumerate}
\item For each open set $U \in \col{T}$, a set $\shf{S}(U)$, called the \emph{stalk at} $U$, 
\item For each pair of open sets $U \subseteq V$, there is a function $\shf{S}(U \subseteq V):\shf{S}(V)\to\shf{S}(U)$, called a \emph{restriction function} (or just a \emph{restriction}), such that
\item For each triple $U \subseteq V \subseteq W$ of open sets, $\shf{S}(U \subseteq W) = \shf{S}(U \subseteq V) \circ \shf{S}(V \subseteq W)$ and
\item $\shf{S}(U \subseteq U)$ is the identity function.
\end{enumerate}
Dually, a \emph{precosheaf $\cshf{C}$ of sets on $(X,\col{T})$} consists of the opposite specification:
\begin{enumerate}
\item For each open set $U \in \col{T}$, a set $\cshf{C}(U)$, called the \emph{costalk at} $U$,
\item For each pair of open sets $U \subseteq V$, there is a function $\cshf{C}(U \subseteq V):\cshf{C}(U)\to\cshf{C}(V)$, called an \emph{extension function} (or just a \emph{extension}), such that
\item For each triple $U \subseteq V \subseteq W$ of open sets, $\shf{C}(U \subseteq W) = \cshf{C}(V \subseteq W) \circ \cshf{C}(U \subseteq V)$ and
\item $\cshf{C}(U \subseteq U)$ is the identity function.
\end{enumerate}

If for every $U \in \col{T}$ there is a pseudometric $d_U$ on the (co)stalk at $U$, and each restriction (or extension) is continuous with respect to the corresponding pseudometrics, we call the entire collection of data a \emph{pre(co)sheaf of pseudometric spaces}.
\end{definition}

As Definition \ref{df:presheaf} makes clear, pre(co)sheaves on a topological space are only sensitive to the poset of open sets, and \emph{not} to the points in those open sets.
In our context, the set of values should be interpreted as the set of values that a collection of variables in a DSEM can take.

\begin{definition}
\label{def:assignment}
    Suppose $\shf{S}$ is a presheaf on a topological space $(X, \col{T})$.
    An \emph{assignment $a$ supported on $\col{U} \subseteq \col{T}$} is an element of the direct product, $\prod_{U \in \col{U}}\shf{S}(U)$.

The direct product is in general \emph{not} the direct sum, since the topology may be infinite!  
For this reason, dually, if $\cshf{C}$ is a precosheaf on $(X,\col{T})$, then a \emph{coassignment supported on $\col{U} \subseteq \col{T}$} is an element of
\begin{equation*}
  \left(\bigsqcup_{U \in\col{U}}\cshf{C}(U) \right).
\end{equation*}

If $\col{U} = \col{T}$, we usually say that the (co)assignment is \emph{global}.
\end{definition}

(Co)assignments may or may not be consistent with their pre(co)sheaf structure.
When they are fully consistent, we highlight this fact by calling them (co)sections.

\begin{definition}
\label{df:section}
A \emph{global section} of a presheaf $\shf{S}$ on a topological space $(X,\col{T})$ is a global assignment $s$ such that for all open $V \subseteq U$ then $\shf{S}(V \subseteq U)\left(s(U)\right) = s(V)$.

Dually, a \emph{global cosection} of a precosheaf $\cshf{C}$ on a topological space is a global coassignment $c$ of the disjoint union under an equivalence,
\begin{equation*}
  \cshf{C}(X) = \left(\bigsqcup_{U \text{ open}}\cshf{C}(U) \right)/\sim,
\end{equation*}
where $\sim$ is the equivalence relation generated by $c_1 \sim c_2$ whenever $c_1 \in \cshf{C}(U_1)$, $c_2 \in \cshf{C}(U_2)$, with $U_1 \subseteq U_2$, and $\left(\cshf{C}(U_1 \subseteq U_2)\right)(c_1) = c_2$.

\emph{Local (co)sections} are defined similarly, but refers to some collection $\col{U}$ of open sets.
\end{definition}

Intuitively, a (co)section corresponds to data that is fully consistent with the hypothesis posed by a (co)sheaf.

The set of global sections of a presheaf on a topological space may be quite different from $\shf{S}(X)$.  It is for this reason that when studying presheaves over topological spaces, an additional \emph{gluing axiom} is included to remove this distinction.
A similar axiom applies for cosheaves.

\begin{definition}
\label{df:sheaf}
Let $\shf{P}$ be a presheaf on the topological space $(X,\col{T})$.  We call $\shf{P}$ a \emph{sheaf on $(X,\col{T})$} if for every open set $U \in \col{T}$ and every collection of open sets $\col{U}\subseteq \col{T}$ with $U = \cup \col{U}$, then $\shf{P}(U)$ is isomorphic to the space of sections over the set of elements $\col{U}$.

Dually, a precosheaf $\cshf{C}$ is a \emph{cosheaf on $(X,\col{T})$} if for every open set $U \in \col{T}$ and every collection of open sets $\col{U}\subseteq \col{T}$ with $U = \cup \col{U}$, then $\cshf{C}(U)$ is isomorphic to the space of cosections over the set of elements $\col{U}$.
\end{definition}

For the time being, we will focus on sheaves.  Cosheaves will reappear in Section \ref{sec:subsystems}.

Given that most assignments are not sections, it is useful to be able to measure how far away an assignment is from being a section.
When we have pseuodmetrics on the stalks, one useful estimate of that distance is the consistency radius.

\begin{definition}
    \label{df:consistency_radius}
    If $\shf{S}$ is a presheaf of pseudometric spaces on a topological space $(X,\col{T})$ and $a$ is a global assignment, the \emph{$p$-norm consistency radius of $a$} is the quantity
    \begin{equation}
    \label{eq:consistency_radius}
c_{\shf{S}}(a) := \left(\sum_{U \in \col{T},}\sum_{V \in \col{T}: V \subseteq U}\left(d_V\left(a(V), \shf{S}(V \subseteq U) a(U) \right)\right)^p\right)^{1/p},
    \end{equation}
    where $p \ge 1$.
\end{definition}

In all of our examples, $p=2$ is used.
A subtle point is that the relative weight of each of the different terms in Equation \eqref{eq:consistency_radius} is implicitly carried by the pseudometrics $d_V$.  
For instance, if $x,y \in \mathbb{R}^n$, a weighted form of the Euclidean pseudometric could be written
\begin{equation*}
    d_V(x,y) = \alpha_V \left(\sum_{k=1}^n |x_k - y_k|^p\right)^{1/p},
\end{equation*}
where $\alpha_V > 0$ is a constant that weighs the importance of the value in the stalk on $V$ in the overall consistency radius.
In some cases, for instance if different units of measure are involved, the correct choice of $\alpha_V$ is clear.
In others, the $\alpha_V$ is a nuisance parameter that needs to be explored by the modeler.

\begin{corollary}
    If $s$ is a global section of a presheaf $\shf{S}$ of pseudometric spaces, then $c_{\shf{S}}(s) = 0$.
\end{corollary}

Consistency radius is stable under perturbations, which means that it can be reliably estimated.

\begin{theorem}\cite[Thm. 1]{robinson2020assignments}
    Consistency radius is a continuous real-valued function of the assignment.    
\end{theorem}

We will often need to consider local assignments as well.  
A natural definition is to define the consistency radius of a local assignment to be the consistency radius of the ``best'' extension of the local assignment to a global one. 

\begin{definition}
\label{df:local_consistency_radius} \cite[Def. 16]{robinson2020assignments}
    If $\shf{S}$ is a presheaf of pseudometric spaces on a topological space $(X,\col{T})$ and $a$ is an assignment supported on $\col{U}\subseteq \col{T}$, then its consistency radius is
    \begin{equation*}
        c_{\shf{S}}(a;\col{U}) := \min \left\{ c_{\shf{S}}(b) : b \in \prod_{U\in\col{T}}\shf{S}(U)\text{ such that } b(U)=a(U)\text{ if }U \in \col{U}\right\}.
    \end{equation*}

    We will use the phrase \emph{minimizing the consistency radius of $a$} as a shorthand for finding the global assignment 
    \begin{equation*}
        b_* := \argmin \left\{ c_{\shf{S}}(b) : b\in \prod_{U\in\col{T}}\shf{S}(U)\text{ such that } b(U)=a(U)\text{ if }U \in \col{U}\right\}.
    \end{equation*}
\end{definition}

As the rest of this chapter shows, minimizing the consistency radius of a given local assignment is the primary tool for sheaf-based inference.

\subsection{The netlist sheaf}
\label{sec:netlist_sheaf}

The key result of this section is that inference for a DSEM corresponds to consistency radius minimization.  In general, it is enabled by Definition \ref{df:netlist_from_dsem} that translates a DSEM into a netlist, and Definition \ref{df:netlist_sheaf} that translates a netlist into a sheaf, in such a way that solutions correspond to global sections (Theorem \ref{thm:dsem_solutions_in_sheaf}).    

\begin{figure}
    \centering
    \includegraphics[width=\linewidth]{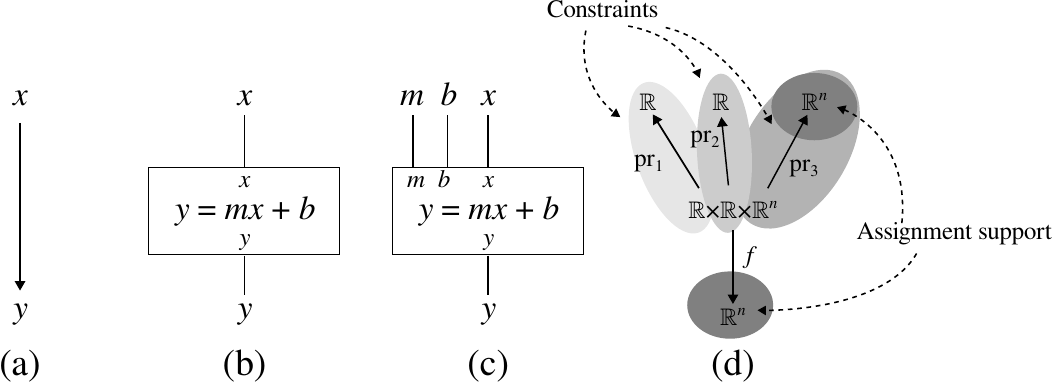}
    \caption{A linear regression problem as (a) a SEM, (b) a netlist with hardcoded coefficients, (c) a netlist with coefficients exposed as inputs, and (d) a sheaf.  To solve the linear regression problem, the partial assignment supported on the darkest shaded region is supplied by the observations, and then the assignment is extended to the remaining stalks.  Finally, the copies of $m$, $b$, and $x$ that should be constrained so that they are identical are shown by the three lighter shadings.}
    \label{fig:linear_regression_netlist}
\end{figure}

In order to motivate the construction, and to explain some of its subtleties, we delay the formal construction (Definition \ref{df:netlist_sheaf}) until after we have discussed two examples.
The first example represents a classic linear regression problem first as a SEM (which is not dynamical), then as a netlist, and finally as a sheaf.  
This progression is summarized in Figure \ref{fig:linear_regression_netlist}.

Before delving into the details, let us
consider the meaning of the arrows shown in Figure \ref{fig:linear_regression_netlist}.
The arrows in each of the frames of Figure \ref{fig:linear_regression_netlist} mean different things.  
In the SEM the arrows have a causal interpretation: the value of $x$ determines that of $y$.  
This interpretation carries over into the netlist, where ports are either inputs or outputs.  

In the sheaf diagram the arrows are functions between the stalks.
Since the stalks represent the set of possible values for each variable,
the functions represented by the arrows will be used to extract data stored on the ports and place them on the nets regardless of whether they are inputs or outputs.
There is no intuitive issue with the outputs. 
An output variable is determined by the data within the part it is attached to. 
However, for an input, the only thing the arrow does is extract the corresponding port's value unmodified.
This seems paradoxical!
The point is that when two parts are connected to each other on a net, 
they both have a claim on what the value of the variable should be.  
If the values correspond to a global section of the sheaf,
this is the assertion that both claims on that variable agree, 
namely the variable produced by the output of one port is the same as the variable that reaches the input port attached to the same net. 

Beginning the example in earnest, suppose that $(x_1,y_1), \dotsc, (x_n,y_n)$ are $n$ points in the plane $\mathbb{R}^2$.
As a modeling choice, we suppose that the $x$ values can be used to predict the $y$ values,
or alternatively that $x$ is an explantory variable and $y$ is a response variable.
If we assert that the model should be linear, we are assuming
\begin{equation*}
    y \approx b + m x,
\end{equation*}
where $b$ and $m$ are parameters to be found.
To express this modeling assumption graphically, we write an arrow $x \to y$, yielding the SEM graph in Figure \ref{fig:linear_regression_netlist}(a).

The netlist for the problem represents the same information as in the SEM.  
As shown in Figure \ref{fig:linear_regression_netlist}(b),
the netlist consists of two variables ($x$ and $y$), and one part (the linear equation that predicts $y$ from $x$).

The prediction process depends on the two parameters $b$ and $m$, which can also be considered as inputs.  
This change results in a netlist with four variables ($x$, $y$, $b$, and $m$) and the same part as before, shown in Figure \ref{fig:linear_regression_netlist}(c).

The sheaf representation of the same system is shown in Figure \ref{fig:linear_regression_netlist}(d).
It is considerably more explicit about variable type information.
The stalk over $m$ and $b$ is $\mathbb{R}$, since each of these parameters takes a real value.
On the other hand, the stalk over $x$ and $y$ is $\mathbb{R}^n$, since they are each a sequence of $n$ real values.
The stalk over the single part is the set of its inputs, namely $\mathbb{R}\times\mathbb{R}\times\mathbb{R}^n$,
corresponding to $m$, $b$, and $x$, respectively.
The restriction maps from the part to the inputs are all projection maps, which select the different inputs.
Explicitly,
\begin{equation*}
    \pr_1\left(m,b,(x_1,\dotsc,x_n)\right) = m,
\end{equation*}
\begin{equation*}
    \pr_2\left(m,b,(x_1,\dotsc,x_n)\right) = b,
\end{equation*}
and
\begin{equation*}
    \pr_3\left(m,b,(x_1,\dotsc,x_n)\right) = (x_1,\dotsc,x_n).
\end{equation*}
The remaining restriction map $f$ shown in Figure \ref{fig:linear_regression_netlist}(d) performs the prediction process, and is given by
\begin{equation}
\label{eq:linear_regression_restriction}
(y_1, \dotsc,y_n) = f\left(m,b,(x_1,\dotsc,x_n)\right) = (mx_1+b, \dotsc, mx_n+b).
\end{equation}
The function $f$ applies the common coefficients ($b$ and $m$) to each of the input values $x_k$ to yield the corresponding output values $y_k$.

The space of global assignments for the sheaf shown in Figure \ref{fig:linear_regression_netlist}(d) is given by the product of all of the stalks.  
This means there are \emph{two copies} of $m$, $b$, and $x$ in the space of global assignments, one for the value of the variable and one as a component of the part.  A typical global assignment $a$ is of the form
\begin{equation}
\label{eq:linear_regression_bundled}
   a := \left(m,b,(x_1, \dotsc, x_n), (y_1, \dotsc, y_n), \left(\widetilde{m},\widetilde{b},(\widetilde{x_1}, \dotsc, \widetilde{x_n})\right)\right),
\end{equation}
where we have listed the four variables first followed by the part.
The consistency radius of this assignment is
\begin{equation}
\label{eq:linear_regression_cr}
c(a) = \left(|\widetilde{m}-m|^p + |\widetilde{b}-b|^p + \sum_{k=1}^n|\widetilde{x_k}-x_k|^p + \sum_{k=1}^n|b+m\widetilde{x_k}-y_k|^p\right)^{1/p}
\end{equation}
for a given $p$.
In what follows, we will take $p=2$, so as to agree with classical linear regression.

The problem of classical linear regression seeks real numbers $m$ and $b$ minimizing the last term in Equation \eqref{eq:linear_regression_cr}.
Therefore, minimizing consistency radius subject to the constraint that each pair of copies of $m$, $b$, and $x$ is equal, and that only $m$ and $b$ are allowed to vary will recover linear regression from the sheaf.  These copies are identified in the lighter shaded regions in Figure \ref{fig:linear_regression_netlist}(d).

\begin{figure}
    \centering
    \includegraphics[width=0.75\linewidth]{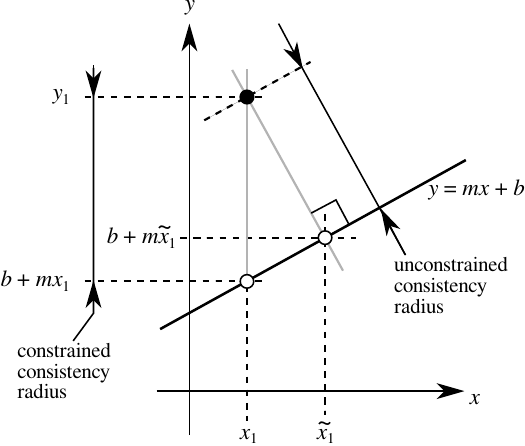} 
    \caption{Geometric meanings of the terms contributing to consistency radius in Equation \ref{eq:linear_regression_cr}.}
    \label{fig:linear_regression_unconstrained}
\end{figure}

To follow the paradigm of consistency radius minimization, we specify a local assignment to the variables $x$ and $y$, and then extend the assignment to a global one.  The support of the local assignment is expressed by the darkest shaded region in Figure \ref{fig:linear_regression_netlist}(d).
Notice that the nets have no higher elements in the partial order shown in Figure \ref{fig:linear_regression_netlist}, so the support of this assignment is $\col{U} = \{\{x\}, \{y\}\}$.
Explicitly, we start with a non-global assignment supported on $\col{U}$,
\begin{equation}
\label{eq:linear_regression_local}
   \left(-,-,(x_1, \dotsc, x_n), (y_1, \dotsc, y_n), -\right),
\end{equation}
where the dashes indicate stalks outside the support of the assignment.
If we seek a global assignment $g^*$ such that
\begin{equation*}
    g^* = \argmin \{c(b) : g(U) = a(U) \text{ for }U \in \col{U}\},
\end{equation*}
this means that we wish to find the entries in the assignment in Equation \eqref{eq:linear_regression_bundled} that are marked with the dashes in Equation \eqref{eq:linear_regression_local}, namely
\begin{equation*}
    \widetilde{m}, \widetilde{b}, m, b, \text{ and } (\widetilde{x_1}, \dotsc, \widetilde{x_n}).
\end{equation*}
Minimizing consistency radius is therefore given by the problem
\begin{equation*}
   \argmin_{\widetilde{m}, \widetilde{b}, m, b, (x_1, \dotsc, x_n)} \left(|\widetilde{m}-m|^2 + |\widetilde{b}-b|^2 + \sum_{k=1}^n|\widetilde{x_k}-x_k|^2 + \sum_{k=1}^n|b+m\widetilde{x_k}-y_k|^2\right)^{1/2}.
\end{equation*}
But since both $\widetilde{m}$ and $m$, and $\widetilde{b}$ and $b$ are being minimized, the consistency radius reduces to
\begin{equation*}
   \argmin_{m, b, (x_1, \dotsc, x_n)} \left(\sum_{k=1}^n|\widetilde{x_k}-x_k|^2 + \sum_{k=1}^n|b+m\widetilde{x_k}-y_k|^2\right)^{1/2}. 
\end{equation*}
This permits the values of the variables $x$ and $y$ to differ from their copies, subject to a penalty.
Instead of least squares regression, this problem is what is usually called \emph{total} least squares; see Figure \ref{fig:linear_regression_unconstrained}.
After minimization, the differences between each of the copies 
\begin{equation*}
    |\widetilde{x_k}-x_k|
\end{equation*}
expresses the uncertainty of their values if the model is to be taken as a given.

To obtain classical least squares regression, we must constrain $\widetilde{x_k}=x_k$ for all $k$.
The global assignment we seek is of the form
\begin{equation*}
    g^* = \left(m,b,(x_1, \dotsc, x_n), (y_1, \dotsc, y_n), \left(m,b,(x_1, \dotsc, x_n)\right)\right),
\end{equation*}
so that the consistency radius minimization problem subject to this constraint becomes
\begin{equation*}
   \argmin_{m, b} \left(\sum_{k=1}^n|b+mx_k-y_k|^2\right)^{1/2}. 
\end{equation*}

Consistency radius minimization unifies several different inference tasks in Figure \ref{fig:linear_regression_netlist}, depending on the support of the initial assignment:
\begin{description} 
    \item[Forward prediction] Choose an assignment supported on $x$, $b$, and $m$, of the form
    \begin{equation*}
        \left(m,b,(x_1, \dotsc, x_n), -, -\right).
    \end{equation*}
    Consistency radius minimization will infer the values for $y$.  Because the above assignment extends to a global section, namely,
    \begin{equation*}
        \left(m,b,(x_1, \dotsc, x_n), (b+mx_1, \dotsc, b+mx_n), \left(m,b,(x_1, \dotsc, x_n)\right)\right),
    \end{equation*}
    consistency radius minimization does not require constraints in this case.
    \item[Backward prediction] Choose an assignment supported on $y$ and $b$, and $m$, of the form
    \begin{equation*}
        \left(m,b,-, (y_1, \dotsc, y_n), -\right).
    \end{equation*}  
    Consistency radius minimization will infer the values for $x$.   If $m \not= 0$, this always results in a global section,
    \begin{equation*}
        \left(m,b,((y_1-b)/m, \dotsc, (y_n-b)/m, (y_1, \dotsc, y_n), \left(m,b,((y_1-b)/m, \dotsc,(y_n-b)/m\right)\right),
    \end{equation*}
    so consistency radius minimization does not require constraints.
    If $m=0$ then the minimizers of consistency radius all have the same consistency radius, and are assignments of the form
    \begin{equation*}
        \left(0,b,(x_1, \dotsc, x_n, (y_1, \dotsc, y_n), \left(0,b,(x_1, \dotsc,x_n)\right)\right).
    \end{equation*}
    Noting that the two copies of the $x$ variable are always identical, applying constraints does not change the result.
    \item[Regression (model fitting)] (Details above, included for completeness here.) Choose an assignment supported on $x$ and $y$, of the form
    \begin{equation*}
        \left(-,-,(x_1, \dotsc, x_n), (y_1, \dotsc, y_n), -\right).
    \end{equation*}
    Consistency radius minimization will infer the values for $b$ and $m$.  As noted above, without constraints consistency radius minimization solves total least squares, while constraints are necessary to recover classical regression.
\end{description}
Hybrid versions of the above problems can also be addressed.

\begin{figure}
    \centering
    \includegraphics[width=0.8\linewidth]{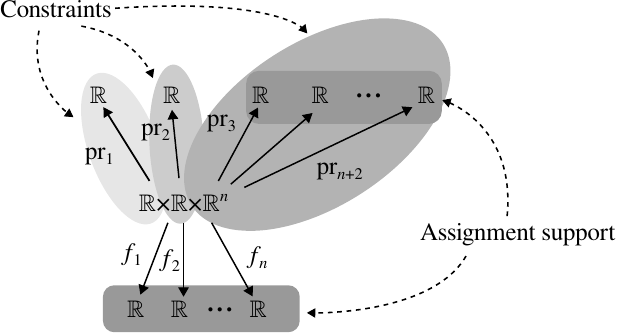}
    \caption{Modification to the sheaf in Figure \ref{fig:linear_regression_netlist}(d) to allow for missing data.}
    \label{fig:linear_regression_sheaf}
\end{figure}

Assignments are populated stalk-wise, so the sheaf in Figure \ref{fig:linear_regression_netlist}(d) explicitly requires that we have access to all of the $n$ data points, since the stalks for $x$ and $y$ are each $\mathbb{R}^n$.  
If there is missing data, a different sheaf construction is possible, in which each separate component of $x$ and $y$ is given its own stalk.  Figure \ref{fig:linear_regression_sheaf} shows the resulting construction.

The $f_k$ restriction maps appearing in Figure \ref{fig:linear_regression_sheaf} are the individual components of the $f$ restriction map in Figure \ref{fig:linear_regression_netlist}(d), namely given Equation \eqref{eq:linear_regression_restriction}, 
\begin{equation*}
    y_k = f_k\left(m,b,(x_1,\dotsc,x_n)\right) = m x_k+b.
\end{equation*}

The set of global assignments for the sheaf in Figure \ref{fig:linear_regression_netlist}(d) is the same as that for the sheaf in Figure \ref{fig:linear_regression_sheaf}, but its components are delineated differently.
A typical global assignment $a$ for the sheaf in Figure \ref{fig:linear_regression_sheaf} is given by 
\begin{equation*}
   a := \left(m,b,x_1, \dotsc, x_n, y_1, \dotsc, y_n,\left(\widetilde{m},\widetilde{b},\widetilde{x_1}, \dotsc, \widetilde{x_n}\right)\right),
\end{equation*}
where the main difference between the above and Equation \eqref{eq:linear_regression_bundled} is in the placement of parentheses.
The consistency radius for a global assignment in both sheaves is given by exactly the same formula.
As in the previous sheaf, we can express the linear regression problem as a consistency radius minimization problem, in which a local assignment supported on the $x_k$ and $y_k$ variables (shown by the darkest shaded regions in Figure \ref{fig:linear_regression_sheaf}) is extended to a global assignment, subject to the constraint that each of the copies of the duplicated variables are identical (shown by the three lighter shaded regions in Figure \ref{fig:linear_regression_sheaf}).
But now, if there is a missing $x_k$ or $y_k$ value, this can simply be excluded from the support of the initial assignment, leaving the specification of the task as a consistency radius minimization unchanged.

\begin{figure}
    \centering
    \includegraphics[width=0.65\linewidth]{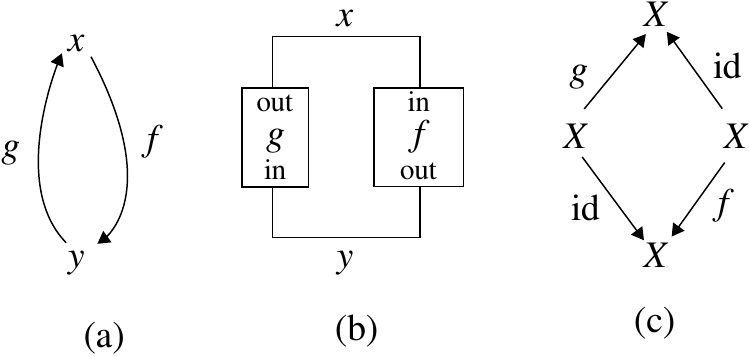}
    \caption{Feedback connections can be handled: (a) a (D)SEM model with feedback, (b) its netlist, (c) its sheaf representation.}
    \label{fig:feedback}
\end{figure}

Feedback connections are easily represented in all of the frameworks under consideration.  Moreover, depending on the set of variables that are permissible, the resulting sheaf will or will not have global sections (Definition \ref{df:section}).

Consider the setting shown in Figure \ref{fig:feedback}:
\begin{description}
    \item[$X=\mathbb{R}$, $f(x)=x$, $g(x)=x$] (Linear SEM) global sections occur whenever the two variables have the same value.
    \item[$X=\mathbb{R}$, $f(x)=-x$, $g(x)=x$] (Linear SEM) the only global section is for both variables to be $0$.
    \item[$X=\mathbb{R}$, $f(x)=1-x$, $g(x)=x$] (Affine, nonlinear SEM) The only global section is for both variables to take the value $1/2$.
    \item[$X=\mathbb{Z}$, $f(x)=1-x$, $g(x)=x$] (Discrete values) No global sections exist.
\end{description}

Feedback will play an important role in defining a sheaf to model autoregressive timeseries in Section \ref{sec:sheaf_autoregression}.

With the preliminary intuition established by the previous two examples, we are now in a position to discuss the general 
translation algorithm.

\begin{definition} 
\label{df:netlist_sheaf}
    If we have a netlist $N$, we build the \emph{netlist sheaf} on the Alexandrov topology of the preorder of its netlist graph of $N$. 
    The stalk on each net is the set of variables for that net.
    The stalk on each part is the product of its input ports.
    The restriction from a part to a net along an input port is the projection function for the corresponding variable set.
    The restriction from a part to a net along an output port is the function that computes the output variable from the set of input variables.
\end{definition}

It is often useful to have individual observations on their own stalks, like we did in Figure \ref{fig:linear_regression_sheaf}.  The following modification to Definition \ref{df:netlist_sheaf} allows for missing data in general.

\begin{definition}
    \label{df:netlist_sheaf_individual}
    Starting with a netlist sheaf as defined in Definition \ref{df:netlist_sheaf}, add an additional element to the preorder of the netlist graph for each observation of each variable.
    These elements are located above their respective variables in the preorder.
    The restriction map from each variable to each observation is the projection that selects the corresponding observation from its parent timeseries.
\end{definition}
 
\begin{theorem}
\label{thm:dsem_solutions_in_sheaf}
Variable values on the netlist correspond bijectively to DSEM solutions and to global sections.
\end{theorem}
\begin{proof}(see also \cite{Robinson_multimodel}[Prop. 6])
    There is a direct correspondence between the values of variables on the nets and the nodes in the DSEM.
    If these are values correspond to a solution, then they directly imply consistency with the restriction maps.
\end{proof}

Moreover, according to \cite[Thm. 1]{robinson2020assignments} there is stability in consistency radius when we perturb away from a consistent set of variables.
This is classical in the case of the linear regression example, because the linear regression coefficients $m$ and $b$ are stable with respect to perturbations in the data variables $x$ and $y$.

\subsection{Sheaves modeling autoregressive timeseries}
\label{sec:sheaf_autoregression}

Autoregressive timeseries are sequences $\dotsc, x_0, x_1, \dotsc$ that obey an equation of the form
\begin{equation*}
    x_n = a_1 x_{n-1} + a_2 x_{n-2} + \dotsb + a_k x_{n-k},
\end{equation*}
for some fixed $a_1, \dotsc, a_k$.
We say that such a sequence is AR($k$) autoregressive.
Autoregressive timeseries can be modeled using the graphical framework being developed in this chapter by the use of feedback connections.

It is easiest to see how the construction of autoregressive timeseries works by starting with a one-step delayed Linear Causal Filter with sliding window size $k$ (which we write as ``LCF($k$)'' for short in diagrams).
Like the linear regression example from the previous section, a variable $x$ is considered an explanatory variable that predicts the values of a response variable $y$. 
This prediction is given by 
\begin{equation*}
    y_n = a_1 x_{n-1} + a_2 x_{n-2} + \dotsb + a_k x_{n-k}
\end{equation*}
where the $a_1, \dotsc a_k$ are constants.

\begin{figure}
    \centering
    \includegraphics[width=0.75\linewidth]{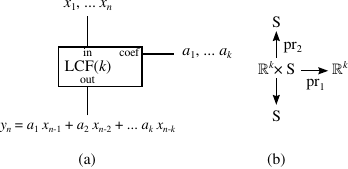}
    \caption{A linear causal filter LCF($k$) with a sliding window size $k$ as (a) netlist wiring hypergraph and (b) netlist sheaf.}
    \label{fig:causal_impulse_sheaf}
\end{figure}

We can realize this equation as a netlist with an input for $x$, an input for $a$, and an output for $y$ shown in Figure \ref{fig:causal_impulse_sheaf}(a).  Using Definition \ref{df:netlist_sheaf}, we obtain the netlist sheaf shown in Figure \ref{fig:causal_impulse_sheaf}(b), where $S$ is the set of infinite sequences of real numbers.

\begin{figure}
    \centering
    \includegraphics[width=\linewidth]{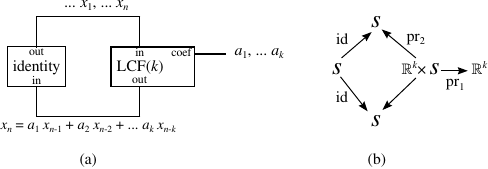}
    \caption{A (a) netlist and (b) a sheaf that encodes an autoregressive timeseries.}
    \label{fig:autoregressive_sheaf}
\end{figure}

To handle autoregressive timeseries, we merely need to consider the pair of equations
\begin{equation*}
\begin{cases}
    y_n = a_1 x_{n-1} + a_2 x_{n-2} + \dotsb + a_k x_{n-k},\\
    x_n = y_n.\\
\end{cases}
\end{equation*}
This is implemented as a netlist with two parts and a feedback connection, as shown in Figure \ref{fig:autoregressive_sheaf}(a), where again $S$ is the set of infinite sequences of real numbers.
The linear causal filter part is the same as before, but the identity part implements the second equation above.
Error terms are not explicitly mentioned, because they are accounted for in the consistency radius calculation (Equation \eqref{eq:consistency_radius}).

The associated netlist sheaf is shown in Figure \ref{fig:autoregressive_sheaf}(b).
Again, consistency radius measures how well the data $x$ fit the model given with coefficients $a$.
Following a theme already present in the linear regression example, there is duplication of data in the sheaf model.
Indeed, the values of $x$ are effectively duplicated in \emph{four} places: the $x$ and $y=x$ variables, and in the two parts.
Once again, if we consider an assignment supported on the two variables (with the same values on each!), minimizing consistency radius will infer the values of the $a$ coefficients. 
Once again, if we run an \emph{unconstrained} optimization, this assumes that some uncertainty is permitted in the values of $x$.

\begin{figure}
    \centering
    \includegraphics[width=\linewidth]{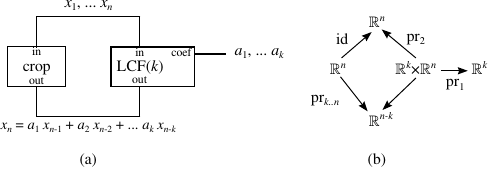}
    \caption{A (a) netlist and (b) a sheaf that encodes an autoregressive timeseries.}
    \label{fig:autoregressive_sheaf_2}
\end{figure}

When the timeseries are finite in length, the equation defining an AR($k$) sequence cannot represent any of the first $k$ time steps.
Therefore, instead of the identity part in Figure \ref{fig:autoregressive_sheaf}, the sheaf for an AR($k$) sequence of length $n$ must crop off the first $k$ components of the vector in the stalk, resulting in a sequence of length $n-k$.
The resulting construction is shown in Figure \ref{fig:autoregressive_sheaf_2}, where we note that a slight abuse of definition occurs in Figure \ref{fig:autoregressive_sheaf_2}(a) because the two outputs are connected to each other.
While this means that the netlist is not valid as such, the sheaf constructed in Figure \ref{fig:autoregressive_sheaf_2}(b) correctly represents an autoregressive sequence.  
Global sections of the sheaf in Figure \ref{fig:autoregressive_sheaf_2}(b) are precisely the AR($k$) sequences of length $n$.

\begin{figure}
    \centering
    \includegraphics[width=\linewidth]{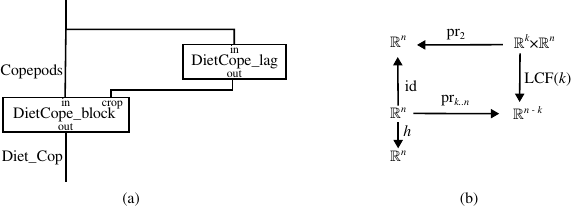}
    \caption{Modification to Figure \ref{fig:schematic_sheafify}(d) to support autoregressive timeseries, shown for the $\mathit{Copepods}$ variable: (a) netlist wiring hypergraph, (b) sheaf diagram.  This modification is performed for each variable in Figure \ref{fig:schematic_sheafify} resulting in Figure \ref{fig:bering_sea_sheaf}.}
    \label{fig:copepods_lagvar}
\end{figure}

Autoregressive sequences can be modeled in the sheaf shown in Figure \ref{fig:schematic_sheafify}(d), our ecological example.
All that is needed is a modification to each variable in the netlist to ensure that each variable is an autoregressive sequence.
Specifically, each of the input variables for each of the parts in the netlist shown in Figure \ref{fig:schematic_sheafify}(b) must be duplicated to represent a lagged copy of the variable, 
and there must be a new part added for each variable to perform the autoregression itself.
As in Figure \ref{fig:autoregressive_sheaf_2}, each original variable gets wired to the input of the corresponding LCF part.
The duplicated (lagged) input on each preexisting part is cropped to be only the most recent samples (since the timeseries is finite), and then that is what is attached to the output port of the LCF part.  
The transformation that is required for the $\mathit{Copepods}$ variable is shown in Figure \ref{fig:copepods_lagvar}.

\section{Sheaf encoding of the Bering Sea}
\label{sec:bering_sheaf}

We now return to the ecological DSEM example introduced in Section \ref{sec:bering_dsem},
and refer the reader to Figure \ref{fig:schematic_sheafify}.
The reader is directed to \cite{netlist_sheaf_github} for the software that generates the sheaf results presented in this section.

The DSEM is shown in Figure \ref{fig:schematic_sheafify}(a), its corresponding netlist wiring hypergraph is shown in Figure \ref{fig:schematic_sheafify}(b), its netlist graph is shown in Figure \ref{fig:schematic_sheafify}(c), and its netlist sheaf is shown in Figure \ref{fig:schematic_sheafify}(d).

\begin{figure}
    \centering
    \includegraphics[width=\linewidth]{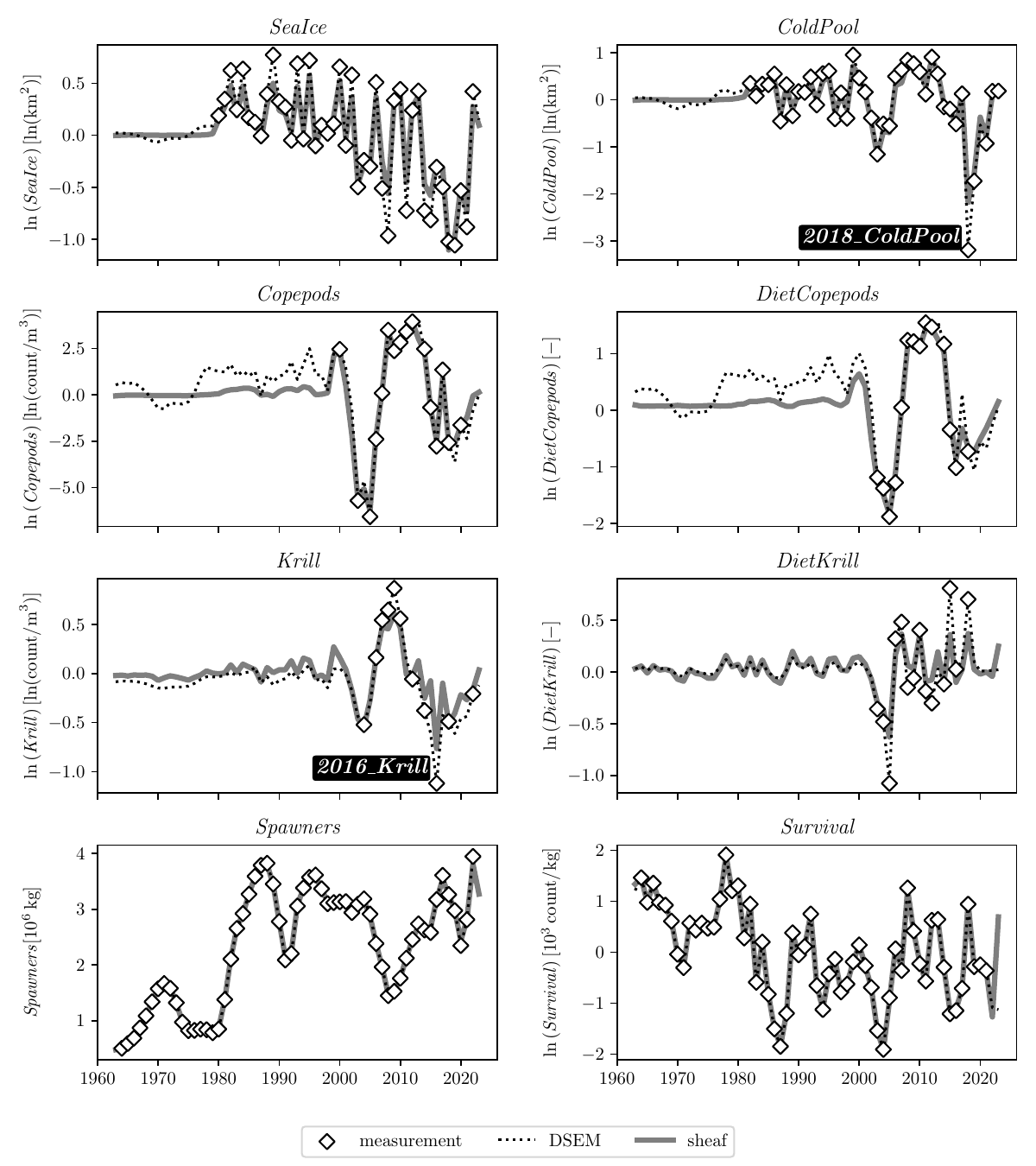}
    \caption{Comparison between the DSEM output and the sheaf with hard-coded path coefficients shown in Figure \ref{fig:schematic_sheafify}(d) and AR($2$) timeseries.  The DSEM was constrained to fit the measurements exactly, whereas the sheaf had no such constraints applied.}
    \label{fig:dsem_sheaf_hardcodedpathcoeffs}
\end{figure}

The netlist sheaf in Figure \ref{fig:schematic_sheafify}(d) does not express the path coefficients as variables, as they are instead ``hard coded'' within each part.
Nevertheless, if the path coefficients are known (for instance, they can be taken from \cite{thorson2024dynamic}), then the sheaf model can be used to predict the values of each of the variables, starting from $\mathit{SeaIce}$ and $\mathit{Spawners}$.
If we apply the modification to the sheaf to require AR($1$) timeseries so that missing data values are interpolated, and use the path coefficients stated in \cite{thorson2024dynamic} (see Table \ref{tab:path_coefficients}), the resulting timeseries are shown in Figure \ref{fig:dsem_sheaf_hardcodedpathcoeffs}.

The DSEM was constrained to fit the measurements exactly, whereas the sheaf had no such constraints applied.
Where the sheaf differs from the measurements, the extent of that difference is a measure of the uncertainty in the value of the variable at the given time. 
This uncertainty is composed of both the measurement and exogenous errors; the sheaf model does not distinguish between the types of error.
Moreover, where there are no measurements available (especially for the earlier measurements), the DSEM reports the expected mean.  The sheaf predictions are typically close to these mean values.
Nevertheless, there is close agreement throughout.
This is not unexpected, because both the sheaf and the DSEM approach are approximations to the same DSEM solution.
There are some differences on the behavior of the earlier inferred data, because many of the observations are missing there.
In these regions, the sheaf tends to yield somewhat less variable predictions than the DSEM (except in the case of the $\mathit{Krill}$ variable). 

As noted earlier, we will compute consistency radius using the Euclidean $p=2$ norm.
Lacking other information, we chose to weight the terms in Equation \eqref{eq:consistency_radius} equally.
The consistency radius of the assignment after minimization is $11.9$.
Since this is not zero, this means that the fit between the data and the model is not perfect.
While the DSEM fits the data for maximum likelihood, the sheaf fits for minimum inconsistency. 
This difference in optimization task results in the observed differences between the sheaf and the DSEM.

\begin{figure}
    \centering
    \includegraphics[width=0.75\linewidth]{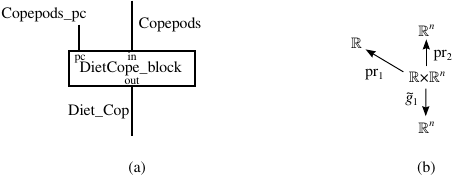}
    \caption{Modification to the netlist to include path coefficients and constants as an input.}
    \label{fig:path_coef_breakout}
\end{figure}

\begin{figure}
    \centering
   \includegraphics[width=\linewidth]{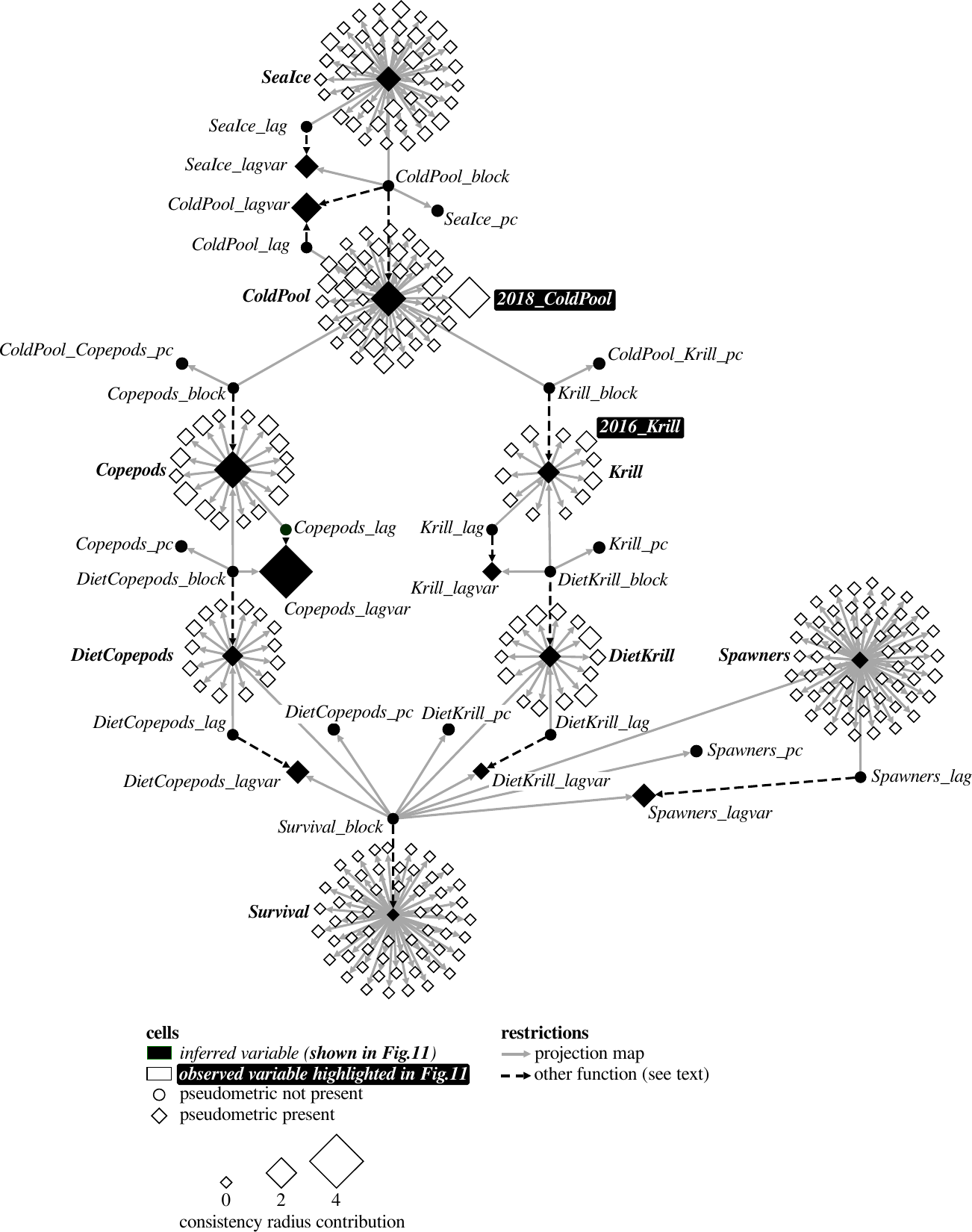}
    \caption{The full sheaf for the DSEM described in Section \ref{sec:bering_dsem}. Its structure reflects the hexagonal backbone shown in the diagrams in Fig.~\ref{fig:schematic_sheafify}.  The black cells represent inferred variables, with the variable names shown in italics.  Variable names that are also bold correspond to variables plotted in Fig.~\ref{fig:dsem_sheaf_hardcodedpathcoeffs}.  White cells represent variables that are observed.  All observed variables except for two are not labeled for clarity.  The two that are labeled have their names in white italics with black backgrounds.  These variables exhibit relatively large contributions to the consistency radius and are highlighted in Fig.~\ref{fig:dsem_sheaf_hardcodedpathcoeffs}.}
    \label{fig:bering_sea_sheaf}
\end{figure}

Taking a cue from Figure \ref{fig:linear_regression_netlist} in the previous section, we can break out path coefficients as separate variables so that they can be adjusted or estimated.
Figure \ref{fig:path_coef_breakout} shows how one of the parts in the netlist shown in Figure \ref{fig:schematic_sheafify}(b) can be modified so that its path coefficients are inputs.  
To handle missing data, we apply Definition \ref{df:netlist_sheaf_individual} to the netlist sheaf, which results in Figure \ref{fig:bering_sea_sheaf}.

Using the sheaf shown in Figure \ref{fig:bering_sea_sheaf}, we can infer the path coefficients and autoregressive coefficients by consistency radius minimization.
Specifically, we construct an assignment supported only on the values of the variables that correspond to observations present in the data.
Then, when we minimize consistency radius, the values of the path coefficients, autoregressive coefficients, \emph{and any missing observations} will be inferred.
The resulting global assignment has a complete timeseries---no missing observations---for each variable as well as path coefficients and autoregressive coefficients.  
Because the approach explained in Section \ref{sec:dsem_background} uses a different strategy for approximating solutions to the problem posed by the DSEM, the inferred path coefficients and missing observations will be somewhat different from those inferred by the sheaf.

There are some differences between the sheaf and the measurement data.
The contributions to consistency radius are not uniformly distributed over the sheaf.
Some of the inconsistency is due to disagreements between the measurements and the DSEM graph model,
and some of the inconsistency is due to the fact that the measurements are not AR($1$) timeseries.  
This is visually apparent in Figure \ref{fig:bering_sea_sheaf}, where it is shown that the two largest contributors to the consistency radius are
\begin{enumerate}
    \item the autoregression cell for $\mathit{Copepods}$ (labeled \emph{Copepods\_lagvar}), and
    \item the year 2018 observations of $\mathit{ColdPool}$ (labeled \emph{2018\_ColdPool}).
\end{enumerate}  

The second of these is easier to interpret.  
We should suspect that the 2018 observation of $\mathit{ColdPool}$ is an outlier (in the $L^2$ sense) from what was expected from the model, and that these differences may have propagated into other parts of the model.
This probably explains why the 2018 observations of $\mathit{Krill}$ and $\mathit{DietKrill}$ are substantially different from the sheaf predictions in Figure \ref{fig:dsem_sheaf_hardcodedpathcoeffs}.

We should interpret the largest contributor to consistency radius as suggesting that the $\mathit{Copepods}$ variable is not well represented by an AR($1$) timeseries.  Notice that the $\mathit{Copepods}$ observations contribute equally to consistency radius, since the small white diamonds encircling the $\mathit{Copepods}$ variable are about the same size. 
This suggests that it is simply that the assumption of Copepods being represented by an AR($1$) timeseries is faulty, rather than any particularly bad observation.

\begin{table}
    \centering
    \begin{tabular}{|c|c||c|c|c|c|c|}
    \hline
    Source & Target & DSEM \cite{thorson2024dynamic}  & \multicolumn{4}{c|}{Sheaf} \\
     &  &  AR($1$) & none & AR($1$) &AR($2$) & AR($10$) \\
    \hline
$\mathit{SeaIce}$ & $\mathit{ColdPool}$ & $0.6$ & $1.68$ & $1.81$ & $1.78$ & $1.74$ \\
$\mathit{ColdPool}$ & $\mathit{Copepods}$ & $1.79$ & $4.45$ & $4.38$ & $4.47$ & $4.17$ \\
$\mathit{ColdPool}$ & $\mathit{Krill}$ & $0.18$ & $0.44$ & $0.38$ & $0.41$ & $0.39$\\
$\mathit{Copepods}$ & $\mathit{DietCopepods}$ & $0.29$ & $0.32$ & $0.35$ & $0.36$ & $0.34$\\
$\mathit{Krill}$ & $\mathit{DietKrill}$ & $0.06$ & $0.52$ & $0.70$ & $0.65$ & $0.56$\\
$\mathit{DietCopepods}$ & $\mathit{Survival}$ & $0.15$ & $-0.50$ & $-0.12$ & $-0.05$ & $-0.32$\\
$\mathit{DietKrill}$ & $\mathit{Survival}$ & $0.13$ & $7.56$ & $5.29$ & $7.19$ & $5.63$\\
$\mathit{Spawners}$ & $\mathit{Survival}$ &  $-0.59$ & $-0.82$ & $-0.65$ & $-0.55$ & $-0.74$\\
\hline
\multicolumn{2}{|c||}{Consistency radius} & $11.9$ & $6.60$ & $9.48$ & $9.03$ & $7.93$ \\
\hline
\multicolumn{2}{|c||}{Runtime (s)} & $2$ & $2848$ & $2637$ & $2679$ & $2907$ \\
    \hline
    \end{tabular}
    \caption{Comparison between path coefficients estimated from the DSEM and the sheaf}
    \label{tab:path_coefficients}
\end{table}

Table \ref{tab:path_coefficients} shows the path coefficients inferred by the DSEM (using maximum likelihood as explained in Section \ref{sec:bering_dsem}) and by the sheaf (using minimum consistency radius).
Table \ref{tab:autoregressive_coefficients} shows the autoregressive coefficients estimated by the sheaf for the AR($1$) and AR($2$) cases.  (The AR($10$) case is not shown for space considerations.)
The DSEM-derived path coefficients were obtained using the assumption of AR($1$) timeseries.
Several different sheaves were constructed with autoregressive sequences of different window sizes.
As a consequence of the construction of consistency radius, minimizing consistency radius infers the following information: (1) missing observations in any variable, (2) all path coefficients, and (3) autoregressive coefficients for each variable.

\begin{table}
    \centering
    \begin{tabular}{|c||c||c|c|}
    \hline
    Variable & AR($1$) & \multicolumn{2}{|c|}{AR($2$)} \\
             &  lag $1$       &       lag $1$               &     lag $2$    \\
    \hline
    $\mathit{ColdPool}$ & 0.582 & 0.480& 0.202\\
$\mathit{SeaIce}$ & 0.361 & 0.287&  0.190\\
$\mathit{Copepods}$ & 0.828 &  1.16 & -0.442\\
$\mathit{Krill}$ & 0.692 & 0.308& 0.411\\
$\mathit{Spawners}$ & 1.01& 1.78& -0.768\\
$\mathit{DietCopepods}$ & 0.886& 1.68 & -0.924\\
$\mathit{DietKrill}$ &0.060 &0.0596 &  0.0445\\
    \hline
    \end{tabular}
    \caption{Autoregressive cofficients estimated by the sheaf for AR($1$) and AR($2$) models.}
    \label{tab:autoregressive_coefficients}
\end{table}

There is broad agreement about the values of the path coefficients between the sheaves with different autoregressive window sizes,
and some agreement between the DSEM and the sheaves.
Since the DSEM does not natively imply a consistency radius, the consistency radius shown for the DSEM is that for the sheaf using AR($1$) timeseries and the hard-coded path coefficients as shown.
Because the consistency radius minimization process on that sheaf cannot adjust the path coefficients---it can only adjust the missing observation values and the autoregressive coefficients---the consistency radius is notably higher in this case.

Some caution in comparing consistency radius across the columns of Table \ref{tab:path_coefficients} is needed.
The number of terms in the consistency radius is the same for each of the sheaves in all but the non-autoregressive case (the fourth column from the left).  This is because the autoregressive coefficients and timeseries are bundled as shown in Figure \ref{fig:autoregressive_sheaf_2}.
Naturally enough, the non-autoregressive sheaf's consistency radius contains no terms pertaining to the autoregressive coefficients, and so is expected to be smaller than the others.
The sheaf column listed as ``none'' means that no autoregressive timeseries assumptions were applied.
Because with no autoregressive assumptions in play, the resulting sheaf diagram is smaller, consequently the consistency radius is smaller.
Interestingly, the consistency radius is smallest for the AR($10$) case, which suggests that more flexibility in the autoregressive coefficients leads to somewhat better prediction accuracy in the measurement data.

Runtimes shown in Table \ref{tab:path_coefficients} are representative when run on an Intel Core Ultra 7 155U at 1.4 GHz with 32 GB RAM.
The process was not memory limited and consumes less than 500 MB RAM.
The sheaf runs roughly $1500$ times slower than the DSEM.
This is because the DSEM solves a sparse linear problem, while the sheaf methodology supports fully nonlinear, non-convex problems.
The sheaf software does not attempt to detect whether the problem is linear, so the consistency radius minimization is always performed as a nonlinear, non-convex optimization problem.  

\section{The topology of subsystems} 
\label{sec:subsystems}

Classically, dynamical systems have been studied using the structure of \emph{invariant sets}.
These are subsets of the space of variable values that are preserved by the action of the dynamical system.
This section shows that invariant sets are one half of a duality pair.  
We can take two different perspectives of a multi-scale dynamical system: invariant sets (which lead to cosheaves) versus subsystems (which lead to sheaves).  

We will establish that a dynamical system induces a \emph{cosheaf of invariant sets}.
The cosheaf of invariant sets breaks the global state of the system into different regimes of behavior, which are parameterized by the open sets of the base space topology.
Conversely, there is also a \emph{sheaf of subsystems} that splits the variables into nested collections that each act independently.

We will formalize the \emph{topology of subsystems} as a finite topological space, by using the Alexandrov topology for a specific preorder (Definition \ref{df:alexandrov}).
Each subsystem corresponds to a preorder element, with composite subsystems ``hooked together'' according to the preorder.
The preorder relation decomposes composite subsystems into their component pieces.   
Intuitively, moving ``up'' in the preorder yields more abstracted ``high-level'' systems.  
This is not entirely compatible with \emph{all} system decompositions in the literature, so caution is advised!    
(The intuition of the presentation here is compatible with Kearney et al. \cite{Kearney_2020}, where the system is modeled as a graph.  
In Kearney et al. \cite{Kearney_2020}, vertices are the loci of state variables, and are ``above'' edges in the preorder constructed in that paper.  Our presentation is also compatible with Steward \cite{Steward_1981}, after transitive closure.)  

\subsection{Dynamical systems}
\label{sec:dyn_cat}

\begin{definition} 
  A \emph{dynamical system} is a continuous bijection $f: S \to S$. The set $S$ in this case is called the set of \emph{states} of the dynamical system.
\end{definition}

It is a classical fact that for a fixed timestep, the solutions to a smooth first order differential equation of the form \eqref{eq:dsem_exact} induce a dynamical system \cite{Strogatz}.
As a consequence, the DSEM, netlist, and sheaf models of the previous sections represent dynamical systems.

\begin{definition}
  For a dynamical system $f: S \to S$, a subset $V \subseteq S$ is called an \emph{invariant set} if
  \begin{equation*}
    f(V) \subseteq V.
  \end{equation*}
\end{definition}

\begin{corollary}
If $V$ is an invariant set of $f: S \to S$, then $f$ restricts to a function $f: V \to V$.
\end{corollary}

\begin{definition}
  Suppose that $A \subseteq B$.  The \emph{inclusion} is the function $i: A \to B$ is a function such that $i(x) = x$ for every $x \in A$.  Notice that $(i|A) \circ i = i$.

  Dually, a \emph{projection} is a function $p: B \to A$ such that $p \circ p = p$ and $p|A = \id_A$.
\end{definition}

\begin{proposition}
  \label{prop:dynamical_warmup}
  Suppose that $U$ and $V$ are two invariant sets for a dynamical system $f: S \to S$ and that $U \subseteq V$.  Then the following diagram
  \begin{equation*}
    \xymatrix{
      U \ar[r]^f \ar[d]_-{i} & U \ar[d]^-{i'}\\
      V \ar[r]_f & V
      }
  \end{equation*}
  commutes, where $i$ and $i'$ are appropriate inclusion maps, which is to say that
  \begin{equation*}
      f \circ i = i' \circ f.
  \end{equation*}
\end{proposition}
\begin{proof}
  Suppose that $x \in U$.  Since $U$ is an invariant set, $f(U) \subseteq U$.  However, since $U \subseteq V$, $x \in V$.  Therefore, $f(x) \subseteq V$ because $V$ is also an invariant set.
\end{proof}

\begin{definition}
  The category $\cat{Dyn}$ of dynamical systems has as its objects dynamical systems.  Each morphism of $\cat{Dyn}$ is a commutative diagram of the form
  \begin{equation*}
    \xymatrix{
      S_1 \ar[r]^{f_1} \ar[d]_{g} & S_1 \ar[d]^{g} \\
      S_2 \ar[r]^{f_2} & S_2 \\
      }
  \end{equation*}
  Composition of morphisms is given by composing the $g$ functions.
\end{definition}

\begin{proposition}
Isomorphisms in $\cat{Dyn}$ are conjugacy classes of dynamical systems. 
\end{proposition}

\subsection{The cosheaf endomorphism of invariant sets}
\label{sec:cosheaf_invariant_sets}

The state space of a dynamical system can be decomposed as the (non-disjoint) union of all its invariant sets.
This collection of invariant sets of a dynamical system is also partially ordered by subset inclusion,
which means that the \emph{collection} of invariant sets can be given an Alexandrov topology.
A cosheaf can be defined to capture the relationship between an invariant set and the invariant sets that contain it.
To this end, the cosheaf identifies duplicate points within these invariant sets with each other.

We begin by observing that the invariance of a collection of subsets with respect to a dynamical system is not necessary to define a cosheaf; it can be constructed generally.

\begin{lemma}
  \label{lem:subset_cosheaf}
  Suppose that $\col{U} \subseteq 2^X$ is an arbitrary collection of subsets of a set $X$.  Consider the inclusion partial order on $\col{U}$, given by $U \le V$ whenever $U \subseteq V$.  Define the following precosheaf $\cshf{C}_{\col{U}}$ on the Alexandrov topology of the inclusion partial order $(\col{U},\le)$:
  \begin{enumerate}
  \item $\cshf{C}_{\col{U}}(U) = U$
  \item $\cshf{C}_{\col{U}}(U \le V) = \cshf{C}_{\col{U}}(U \subseteq V): U \to V$ via the inclusion map.
  \end{enumerate}
  Then $\cshf{C}_{\col{U}}$ is a cosheaf of sets on the Alexandrov topology of the inclusion partial order $(\col{U},\le)$.
\end{lemma}
\begin{proof}
  Suppose that $V \in \col{U}$, and that $\col{V} \subseteq \col{U}$ is a collection of subsets with $V = \cup \col{V}$.  We need to establish that the space of global cosections on $\col{V}$ is identical to $\cshf{C}_{\col{U}}(V) = V$.  The space of global cosections on $\col{V}$ is
  \begin{equation*}
    \left(\bigsqcup_{W \in \col{V}} \cshf{C}_{\col{U}}(W) \right)/\sim = \left(\bigsqcup_{W \in \col{V}} W \right)/\sim = \bigcup_{W \in \col{V}} W = \cup \col{V}= V,
  \end{equation*}
  since the equivalence $\sim$ identifies points that agree on overlaps.
\end{proof}

The above cosheaf construction is functorial, which means that it is compatible with transformations of the underlying sets.  
In order to establish functoriality, we need to formalize these transformations by defining the class of morphisms for sheaves and cosheaves.

\begin{definition}
\label{df:morphisms}
Suppose that $\shf{R}$ is a sheaf on $(X,\col{T}_X)$, $\shf{S}$ is a sheaf on $(Y,\col{T}_Y)$, and that $f:(X,\col{T}_X) \to (Y,\col{T}_Y)$ is a continuous function.
A \emph{sheaf morphism $m : \shf{R} \to \shf{S}$} is a collection of maps $m_U : \shf{R}(f^{-1}(U))\to \shf{S}(U)$ for each $U \in \col{T}_Y$ such that the following diagram commutes whenever $U,V \in \col{T}_Y$ and $U \subseteq V$,
\begin{equation*}
    \xymatrix{
    \shf{R}(f^{-1}(V)) \ar[r]^-{m_V} \ar[d]_{\shf{R}(f^{-1}(U) \subseteq f^{-1}(V))} & \shf{S}(V) \ar[d]^{\shf{S}(U \subseteq V)} \\
    \shf{R}(f^{-1}(U)) \ar[r]_-{m_U} & \shf{S}(U) \\    
    }
\end{equation*}
Dually, if $\cshf{R}$ is a cosheaf on $(X,\col{T}_X)$, and $\cshf{S}$ is a cosheaf on $(Y,\col{T}_Y)$, a \emph{cosheaf morphism $m : \cshf{R} \to \cshf{S}$} is a collection of maps $m_U : \cshf{R}(f^{-1}(U)) \to \cshf{S}(U)$ such that the following diagram commutes whenever $U,V \in \col{T}_Y$ and $U \subseteq V$,
\begin{equation*}
    \xymatrix{
    \cshf{R}(f^{-1}(V)) \ar[r]^-{m_V} & \cshf{S}(V) \\
    \cshf{R}(f^{-1}(U)) \ar[r]_-{m_U} \ar[u]^{\cshf{R}(f^{-1}(U) \subseteq f^{-1}(V))} & \cshf{S}(U) \ar[u]_{\cshf{S}(U \subseteq V)}\\    
    }
\end{equation*}
\end{definition}

With the definition of morphisms in hand, we can now establish that the cosheaf construction in Lemma \ref{lem:subset_cosheaf} is functorial.

\begin{lemma}
  \label{lem:top_cosheaf}
  There is a functor $\cat{Top} \to \cat{CoShv}$ that takes a topological space $(X,\col{T})$ to a cosheaf $\cshf{C}_{(X,\col{T})}$ of sets on $(X,\col{T})$ via $\cshf{C}_{(X,\col{T})}(U) := U$ and $\cshf{C}_{(X,\col{T})}(U \subseteq V)$ is the inclusion $U \hookrightarrow V$.
\end{lemma}
\begin{proof}
  First, we observe that Lemma \ref{lem:subset_cosheaf} establishes that $\cshf{C}_{(X,\col{T})}$ is a well-defined cosheaf on $(X,\col{T})$.

  Suppose that $f: (X,\col{T}_X) \to (Y,\col{T}_Y)$ is a continuous map.  This lifts to a cosheaf morphism $F : \cshf{C}_{(X,\col{T}_X)} \to \cshf{C}_{(Y,\col{T}_Y)}$.   Suppose that $U \subseteq V$ are two open sets in $Y$.  Then we have that $f^{-1}(U) \subseteq f^{-1}(V)$ are two open sets in $X$.  Therefore, the following diagram commutes
  \begin{equation*}
    \xymatrix{
      \cshf{C}_{(X,\col{T}_X)}(f^{-1}(U)) = f^{-1}(U) \ar[d]_{\cshf{C}_{(X,\col{T}_X)}(f^{-1}(U) \subseteq f^{-1}(V))} \ar[rr]^-{F_{U}:=f|U} && \cshf{C}_{(Y,\col{T})}(U)= U \ar[d]^{\cshf{C}_{(Y,\col{T}_Y)}(U \subseteq V)} \\
      \cshf{C}_{(X,\col{T}_X)}(f^{-1}(V))=f^{-1}(V) \ar[rr]_-{F_V:=f|V} && \cshf{C}_{(Y,\col{T}_Y)}(V)=V
      }
  \end{equation*}
  which establishes definitions for the component maps of $F$, and therefore that $F$ is a cosheaf morphism.

  Now suppose that we have two continuous maps $f: (X,\col{T}_X) \to (Y,\col{T}_Y)$ and $g:(Y,\col{T}_Y) \to (Z,\col{T}_Z)$.  We must show that the corresponding composition of cosheaf morphisms $G \circ F$ is the equal to the one induced by $(g \circ f)$. This follows immediately because the components maps of the cosheaf morphism $G \circ F$ are simply restrictions of the composition $(g \circ f)$.  
\end{proof}

Suppose that $f: S \to S$ is a dynamical system.  The invariant sets of $f$ are indeed a collection of subsets, which are partially ordered by inclusion.  Therefore, Lemma \ref{lem:subset_cosheaf} establishes that there is a well-defined cosheaf $\cshf{S}$ of invariant sets of $f$.

\begin{proposition}
  \label{prop:invariant_set_endomorphism}
   A dynamical system $f: S \to S$ induces an morphism $m:\cshf{S} \to \cshf{S}$ on the cosheaf of invariant sets, and for which the induced map on global cosections is $m_S = f$.
\end{proposition}
\begin{proof}
  Suppose that $U$ is an invariant set of $f$.  Let $m_U : U \to U$ be the restriction of $f$ to $U$.  If $U \subseteq V$ are two invariant sets, then Proposition \ref{prop:dynamical_warmup} implies that
  \begin{equation*}
    \xymatrix{
      U \ar[r]^-{m_U=f} \ar[d]_{i} & U \ar[d]^{i}\\
      V \ar[r]_-{m_V=f} & V
      }
  \end{equation*}
  commutes, where $i$ is the inclusion map.  It is immediate that this is exactly the condition that the $m$ maps are the components of a cosheaf morphism.  Moreover, since $S$ is itself an invariant set, the proof is complete.
\end{proof}

\subsection{Subsystem decomposition sheaf}
\label{sec:subsystem_sheaf}

Rather than carving up the state space into different regimes of behavior, we can instead carve it into non-interacting collections of variables.
In this way, we arrive at the \emph{subsystem sheaf} instead of the invariant set cosheaf.
The global sections combine variables together into vectors, whereas global cosections paste subsets of values together.

Dualizing the condition for an invariant set yields the condition for a subsystem.  Suppose that $f: S \to S$ is a bijection and that $U \subseteq S$ is an invariant set for $f$.  If $i: U \to S$ is the inclusion map, then the diagram at left below commutes:
\begin{equation*}
  \xymatrix{
    S \ar[r]^f & S & S \ar[d]_p \ar[r]^f & S \ar[d]^p\\
    U \ar[u]^{i} \ar[r]_{f|U} & U \ar[u]_{i} & B \ar[r]_{g} & B\\
    }
\end{equation*}
Dually, the diagram at right above captures the situation where $B$ is a \emph{subsystem} of $f$.  

\begin{definition}
\label{def:subsystem}
  If $f: S \to S$ is a dynamical system,
  a \emph{subsystem} is a pair $(g,p)$ consisting of a dynamical system $g : B \to B$ and a surjection $p : S \to B$ such that $p \circ f = g \circ p$.  We will call $p$ the \emph{subsystem projection}.
When $p$ is clear from context, we will often say \emph{$g$ is a subsystem of $f$}.
\end{definition}

We can think of the function $g$ as a dynamical system in its own right.

The idea of a subsystem is neatly compatible with the DSEM construction.  
As will be shown later in Corollary \ref{cor:acyclic_dsem_subsystems}, when the DSEM graph is acyclic, the subsystems can be ``read off'' directly.
For the moment, a few examples will build the necessary intuition.

\begin{example}
\label{eg:dsem_subsystems_1}
  Consider the DSEM with two variables $A$ and $B$, given by the graph with one edge $A \to B$.  
  The variable $A$ is a subsystem on its own, whereas $B$ cannot be a subsystem on its own because its value cannot be predicted from $B$ alone.  As a result, there are two nested subsystems: $\{A\}$ and $\{A\to B\}$.

  To see this explicitly, suppose that the values of $A$ are given by the timeseries $\{a_n\}$ and the values of $B$ are given by the timeseries $\{b_n\}$, with the prediction of $B$ from $A$ given by the formula
  \begin{equation*}
      b_{n+1} = \beta(a_n, a_{n-1}, \dotsc).
  \end{equation*}
  The dynamical system implied by this DSEM is represented by shifting the timeseries by one timestep.
  Specifically, the dynamical system is given by the function $f: A \times B \to A \times B$ given by
  \begin{equation*}
  \begin{aligned}
      f(\dotsc,&a_n,a_{n-1},\dotsc, \dotsc,b_n,b_{n-1},\dotsc) \\&= (\dotsc,a_{n+1},a_n,\dotsc,\dotsc, \beta(a_n, a_{n-1}, \dotsc),\beta(a_{n-1}, a_{n-2}, \dotsc),\dotsc).
  \end{aligned}
  \end{equation*}
  Because of this formula, it should be clear that $\{B\}$ cannot be a subsystem because the values of the $\{b_n\}$ timeseries depend on the values of $\{a_n\}$.  Under a projection that removes the $\{a_n\}$ from the domain, the values of $\{b_n\}$ cannot be determined.
  
  The subsystem $\{A\}$ arises using the subsystem projection $p: A\times B \to A$, namely
  \begin{equation*}
   p(\dotsc,a_n,a_{n-1},\dotsc, \dotsc,b_n,b_{n-1},\dotsc) = (\dotsc,a_{n+1},a_n,\dotsc).
  \end{equation*}
  The subsystem dynamical map $g: A \to A$ is simply
  \begin{equation*}
   g(\dotsc,a_n,a_{n-1},\dotsc) = (\dotsc,a_{n+1},a_n,\dotsc).
  \end{equation*}
  Verification that $(g,p)$ is a subsystem is then simply a calculation,
  \begin{equation*}
      \begin{aligned}
          (p \circ f) &(\dotsc,a_n,a_{n-1},\dotsc, \dotsc,b_n,b_{n-1},\dotsc) \\ &= p(\dotsc,a_{n+1},a_n,\dotsc,\dotsc, \beta(a_n, a_{n-1}, \dotsc),\beta(a_{n-1}, a_{n-2}, \dotsc),\dotsc)\\
          &=(\dotsc,a_{n+1},a_n,\dotsc)\\
          &= g(\dotsc,a_n,a_{n-1},\dotsc)\\
          &= (g \circ p) (\dotsc,a_n,a_{n-1},\dotsc, \dotsc,b_n,b_{n-1},\dotsc).\\
      \end{aligned}
  \end{equation*}
\end{example}

\begin{example}
\label{eg:dsem_subsystems_2}
  \begin{equation*}
      \xymatrix{
      & B \\
      A \ar[ur] \ar[dr] & \\
      & C \\
      }
  \end{equation*}
  Following the logic of Example \ref{eg:dsem_subsystems_1}, the subsystems are $\{A\}$, $\{A\to B\}$, $\{A \to C\}$, and the original system. 
\end{example}

\begin{example}
  \label{eg:dsem_subsystems_3}
  Consider the DSEM with three variables $A$, $B$, and $C$ given by the graph
  \begin{equation*}
      \xymatrix{
      A \ar[dr] & \\
      & C \\
      B \ar[ur] & \\
      }
  \end{equation*}
  Following the logic of Example \ref{eg:dsem_subsystems_1}, the subsystems are $\{A\}$, $\{B\}$, and the original system.
  Notice that $\{C\}$ cannot be a subsystem on its own because its values are determined by both $A$ and $B$.
\end{example}

When a dynamical system is described by a DSEM with feedback, there are often fewer subsystems because the values of the variables cannot be determined in isolation.

\begin{example}
  \label{eg:logic_oscillator}
  Consider the DSEM on variables $A$ and $B$ given by the graph
  \begin{equation*}
      \xymatrix{
      A \ar@/^/[r] & B \ar@/^/[l] \\
      }
  \end{equation*}
  (See also Figure \ref{fig:feedback} for the sheaf model.)
  In this case, the only subsystem is the entire system, because the values of $A$ cannot be determined without knowing $B$, and conversely the values of $B$ cannot be determined without knowing $A$.
\end{example}

Linear systems are special because invariant sets and subsystems reduce to the same thing, as the next example shows.

\begin{example}
\label{eg:linear_subsystems}
Let $V$ be a finite dimensional vector space and $f: V \to V$ be a linear isomorphism.
If we use the usual Euclidean norm on $V$, $f$ is continuous, so it is also a dynamical system.
Subsystems and invariant subspaces of $f$ are in bijective correspondence.

To see this, suppose that $v \in V$ is an eigenvector for $f$, namely
\begin{equation*}
    f(v) = \lambda v
\end{equation*}
for some $\lambda$.
Then the subspace spanned by $v$ is an invariant set.
Conversely, every invariant set of $f$ is a linear subspace, spanned by a set of eigenvectors (possibly with complex eigenvalues).

Since $V$ was assumed to be finite dimensional, every subspace $W \subseteq V$ also has an associated orthogonal projection $\pr_W : V \to W$.
If $W$ is an invariant set for $f$, then $(f|W,\pr_W)$ is a subsystem.
To see this, suppose that $v \in V$, which can be written as the decomposition $u + w$, where $w \in W$ and $\pr_W(u) = 0$.
Because $f$ is a linear isomorphism, the assumption on $u$ means that $\pr_W(f(u)) = 0$.
All that remains is to verify that the definition of subsystem holds,
\begin{equation*}
    \begin{aligned}
   (\pr_W \circ f)(v) & = \pr_W \left(f(u+w)\right)\\
        & = \pr_W \left(f(u) + f(w)\right)\\
        & = \pr_W \left(f(u)\right) + f(w)\\
        &= f(w)\\
        &= (f|W)\left(w\right) \\
        &= (f|W)\left(\pr_W(u+w)\right) \\
        & = (f|W \circ \pr_W)(v).\\
    \end{aligned}
\end{equation*}
\end{example}

\begin{lemma}
  \label{lem:subsystem_preorder}
  The relation ``is a subsystem of'' is a preorder, or in other words a reflexive, transitive relation.
\end{lemma}
\begin{proof}
  Suppose that $f:S \to S$ is a dynamical system.  Reflexivity follows immediately by taking $(f,\id_S)$ as a subsystem.
  For transitivity, suppose that $(g_2,p_2)$ is a subsystem of $f$, and that $(g_1,p_1)$ is a subsystem of $g_2$.  That is, we have the commutative diagram
  \begin{equation*}
  \xymatrix{
    S \ar[rr]^f \ar[d]^{p_2} \ar@/_1pc/[dd]_{p_1 \circ p_2} && S \ar[d]_{p_2} \ar@/^1pc/[dd]^{p_1\circ p_2} \\
    B_2 \ar[d]^{p_1} \ar[rr]^{g_2} && B_2 \ar[d]_{p_1}\\
    B_1 \ar[rr]_{g_1} && B_1
    }
  \end{equation*}
  so that $(g_1,(p_1 \circ p_2))$ is a subsystem of $f$.
\end{proof}

Intuitively, the preorder specifies how data can flow from one subsystem to the next.  If $(g_1,p_1)$ is a subsystem of $(g_2,p_2)$, then each variable in $(g_2,p_2)$ is also a variable of $(g_1,p_1)$.  As a result, the state of $g_1$ can influence the state of $g_2$.

\begin{example}
  Consider the dynamical system $f: \mathbb{Z}^3 \to \mathbb{Z}^3$ given by
  \begin{equation*}
    f(x,y,z) := ((1-x),y(1-x)+zx,z(1-x)+yx).
  \end{equation*}
  This has a nontrivial subsystem $\pr_{1} : \mathbb{Z}^3 \to \mathbb{Z}$, since the map
  \begin{equation*}
    g(x) := 1-x
  \end{equation*}
  makes the following diagram commute
  \begin{equation*}
    \xymatrix{
      \mathbb{Z}^3 \ar[r]^f \ar[d]_{\pr_1} & \mathbb{Z}^3 \ar[d]^{\pr_1} \\
      \mathbb{Z} \ar[r]_g & \mathbb{Z}\\
      }
  \end{equation*}
  In this case, the $x$ variable in the subsystem acts as an input to the overall system, even though its behavior is isolated from the rest of the system.
\end{example}

It is not necessarily the case that subsystems are invariant sets.

\begin{example}
  Consider the dynamical system $f: \mathbb{R}^2 \to \mathbb{R}^2$, given by $f(x,y) := (x,y+1)$.  Consider the subset $B = \{(x,0) : x \in \mathbb{R}\}$.  This set yields a subsystem, since the following diagram commutes
  \begin{equation*}
    \xymatrix{
      \mathbb{R}^2 \ar[d]_p \ar[r]^f & \mathbb{R}^2 \ar[d]^p \\
      B \ar[r]_{\id} & B
      }
  \end{equation*}
  where $p(x,y) = (x,0)$, even though the set $B$ is not an invariant set.
\end{example}

However, conversely, invariant sets of subsystems do determine invariant sets of their parent system.

\begin{lemma}
  Suppose that $f: S\to S$ is a dynamical system with $g: B \to B$ is a subsystem with subsystem projection $p: S \to B$.  If $V \subseteq B$ is an invariant set of $g$, then $p^{-1}(V)$ is an invariant set of $f$.
\end{lemma}
\begin{proof}
  The hypotheses posit a commutative diagram of the form
  \begin{equation*}
    \xymatrix{
      S \ar[r]^{f} \ar[d]_{p} & S \ar[d]^{p} \\
      B \ar[r]_g & B\\
      }
  \end{equation*}

  Suppose that $x \in p^{-1}(V) \subseteq S$.  We have that $p(f(x)) = g(p(x))$ via the commutative diagram above.  Noting that $p(x) \in V$ by construction, and that $V$ is an invariant set of $g$, this means that $g(p(x)) \in V$.  Thus, $p(f(x)) \in V$, so $f(x) \in p^{-1}(V)$, which establishes that $p^{-1}(V)$ is an invariant set of $f$. 
\end{proof}

\begin{lemma}
  \label{lem:subsystem_descend}
  Suppose that $f: S \to S$ is a dynamical system and that $Y \subseteq S$ is an invariant set for $f$.  If $g: B \to B$ is a subsystem of $f$ with subsystem projection $p$, then $g$ is also a subsystem of $f|Y$.
\end{lemma}
\begin{proof}
  Suppose that $i: Y \to S$ is the inclusion map.  The hypotheses state that the diagram of solid arrows below commutes:
  \begin{equation*}
    \xymatrix{
      Y \ar@{-->}[ddr] \ar[rr]^{(f|Y)} \ar[dr]^i & & Y \ar@{-->}[ddr] \ar[dr]^i &\\
      & S \ar[rr]^{f} \ar[d]^p && S \ar[d]^p \\
      & B \ar[rr]_g && B\\
      }
  \end{equation*}
  The conclusion follows by completing the diagram's dashed arrows with the composition $p \circ i$ as the subsystem projection for $g$ as a subsystem of $f|Y$.
\end{proof}

A related statement to Lemma \ref{lem:subsystem_descend} could consider the conditions under which a subsystem of an invariant set lifts to a subsystem of the entire system.  Diagrammatically, this consists of a situation where the subsystem projections defined by the dashed arrows in the diagram below could be constructed:
\begin{equation*}
  \xymatrix{
    Y \ar[ddr] \ar[rr]^{(f|Y)} \ar[dr]^i & & Y \ar[ddr] \ar[dr]^i &\\
    & S \ar[rr]^{f} \ar@{-->}[d] && S \ar@{-->}[d] \\
    & B \ar[rr]_g && B\\
  }
\end{equation*}

Therefore, when studying a dynamical system, one will often encounter problems of the following form.
\begin{question}
  When do lifts to the dashed arrows in the diagram above exist?
\end{question}

Answers to this question relate closely to the expected behavior of systems when they are rewritten with new variables.
This routinely happens with compiled software, as the next example shows.

\begin{example}
  \label{eg:compiler_commuting}
  Suppose that $X$ represents the state space of a computer, perhaps a Turing machine.  The design of the computer and physical laws yield a dynamical system $f: X \to X$.  For this example, $f$ is not bijective.  

  The way that the computer is used is that the user loads an executable and then runs it.  The initial state of the executable is a point within a subset $U \subseteq X$.  The user does not have control over the entire state of the machine, but rather can constrain it to a smaller portion of the state space.  It makes sense to require that $U$ is an invariant set, which means that not only the initial state is included, but all possible future states as well.  Therefore, the execution of the executable is completely determined by the commutative diagram
  \begin{equation*}
    \xymatrix{
      U \ar[r]^{f|U} \ar[d]_{\ell} & U \ar[d]^{\ell}\\
      X \ar[r]_f & X\\
      }
  \end{equation*}
  
  As an example in PDP-11 assembly, we could have
  \begin{equation*}
    U = \{\texttt{PC} \in \{0,1\}, \text{memory} = \{\texttt{0} : \texttt{ADD R1,R2}, \texttt{1} : \texttt{HALT}\} \},
  \end{equation*}
  where all values of the unspecified parts of the machine state (other registers, the rest the memory) are included in $U$.
  If the program counter {\tt PC} is initialized to $0$, the program will execute the instructions at $0$ and $1$, and then will halt.  Evidently, if $\texttt{PC}=1$, then the program halts immediately.  No modifications to memory can occur given an initialization with $U$, and {\tt PC} cannot be moved outside of those two instructions.  This ensures that $f(U) \subseteq U$ is indeed an invariant set.

  We might instead imagine that the executable specified by $U$ was the result of a compiled, high-level program.  Such a program would necessarily be of the form $g: Y \to Y$, where $Y$ holds the values of the two registers {\tt R1} and {\tt R2}.  For a PDP-11, this means $Y = (\{0,1\}^{16})^2$, and
  \begin{equation*}
    g(x,y) := (x,x+y),
  \end{equation*}
  which is to say that {\tt R1} is unchanged by the program, and {\tt R2} takes the sum of {\tt R1} and {\tt R2}.

  The compilation process essentially ensures that we have the following commutative diagram
  \begin{equation*}
    \xymatrix{
      U \ar[r]^{f|U} \ar[d]_{q} & U \ar[d]^{q}\\
      Y \ar[r]_g & Y\\
      }
  \end{equation*}
  where the $q$ maps select the two registers {\tt R1} and {\tt R2} from the entirety of the machine state.

  Notice that we may write $q = p \circ \ell$, where $\ell$ is the inclusion of $U \hookrightarrow X$, and $p$ still selects the two registers {\tt R1} and {\tt R2} from the entirety of the machine state.  Since the machine state is very large in comparison to $U$, the following diagram does \emph{not} commute:
  \begin{equation*}
    \xymatrix{
      X \ar[r]^{f} \ar[d]_{p} & U \ar[d]^{p}\\
      Y \ar[r]_g & Y\\
      }
  \end{equation*}
  Values of $X$ for which the commutativity fails egregiously are instances of \emph{weird machine states} \cite{Dullien_weird}.

  However, when the operating system loads an executable, there are conventions about initialization.  This helps to avoid weird machine states.  We can formalize this idea by way of an initialization function $i : Y \to U$ that is a right inverse to $q$, namely $q \circ i = (p \circ \ell) \circ i = \id_Y$.  This means that we have the following commutative diagrams
  \begin{equation*}
    \xymatrix{
      U \ar[r]^{f|U}  & U \ar[d]^{q}& X \ar[r]^f & X \ar[d]^{p}\\
      Y \ar[r]_g \ar[u]^i & Y&Y\ar[u]^{\ell \circ i} \ar[r]_g & Y\\
      }
  \end{equation*}
  
  For instance, in the example PDP-11 program, we could use
  \begin{equation*}
    \begin{aligned}
      i(x,y) := \{\texttt{PC} &=0, \\
      \texttt{R1}&=x,\\
      \texttt{R2}&=y,\\
      \texttt{R[3-6]}&=0,\\
      \text{memory}&=\{\texttt{0} : \texttt{ADD R1,R2}, \texttt{1} : \texttt{HALT}, \texttt{[2-]} : 0\} \},\\
      \end{aligned}
  \end{equation*}
  Notice that since $i$ does not have the ability to change the program counter {\tt PC}, the following diagram does \emph{not} commute
  \begin{equation*}
    \xymatrix{
    U \ar[r]^{f|U} & U \\
    Y \ar[u]^i \ar[r]_g & Y\ar[u]_i\\
    }
  \end{equation*}
\end{example}

Inspired by Example \ref{eg:compiler_commuting}, suppose that we have a commutative diagram
\begin{equation*}
  \xymatrix{
    X \ar[r]^f & X\ar[d]^p\\
    Y \ar[u]^i \ar[r]_g & Y\\
    }
\end{equation*}
where $i$ is injective, $p$ is surjective, and $f$, $g$ are bijective.

This leads to another question that is often of interest when studying system behaviors.
\begin{question}
  Under what conditions does
\begin{equation*}
    \xymatrix{
    X \ar[r]^f \ar[d]_p & X\ar[d]^p\\
    Y\ar[r]_g & Y\\
    } 
\end{equation*}
commute?  Clearly if $g$ is bijective, then a sufficient condition is that $p = g^{-1} \circ p \circ f$.  It is probably the case that $p \circ i = \id_Y$ in most applications, but it is unlikely to be the case that $i \circ p = \id_X$.
\end{question}

\begin{lemma}
  \label{lem:subsystem_meet_semi}
  The subsystem preorder is a meet-semilattice.  That is, if we have two subsystems $f_i: S_i \to S_i$ for $i=1,2$ of a dynamical system $f: S\to S$, there is a common subsystem $f_3: S_3 \to S_3$ of both of them (which might be trivial) that satisfies the following universal property.  If $f_4: S_4 \to S_4$ is another common subsystem of $f_1$ and $f_2$, then $f_4$ is a subsystem of $f_3$.
\end{lemma}
\begin{proof}
  We start with two subsystems of a common dynamical system $f: S \to S$, so that we have a commutative diagram
  \begin{equation*}
    \xymatrix{
      S_1 \ar[r]^{f_1} & S_1 \\
      S \ar[u]^{p_1}\ar[d]_{p_2} \ar[r]^f & S \ar[u]_{p_1}\ar[d]^{p_2}\\
      S_2 \ar[r]_{f_2} & S_2
      }
  \end{equation*}
  We want to construct a subsystem of all three of these $f_3: S_3 \to S_3$, that is as large as possible.  Realize that what is needed to satisfy the universal property is a definition for the dashed arrows in
  \begin{equation*}
    \xymatrix{
      S \ar[d]_{p_2} \ar[r]^{p_1} & S_1 \ar@{-->}[d]^{p_3'} \\
      S_2 \ar@{-->}[r]_{p_3''} & S_3\\
      }
  \end{equation*}
  such that this diagram is a colimit.

  Since each of the $S_i$ are sets, there is a standard colimit construction, namely $S_3 = (S_1 \sqcup S_2) / \sim$ where $x \sim y$ if $x \in S_1$, $y \in S_2$ such that there is a $z \in S$ with $p_1(z) = x$ and $p_2(z)=y$.  The colimit condition implies that when we apply this construction twice, there is a unique $f_3$ completing the diagram below
  \begin{equation*}
        \xymatrix{
      S \ar[d]_{p_2} \ar[r]^{p_1} & S_1 \ar[d]^{p_3'} \ar[dr]^{f_1} \\
      S_2 \ar[r]_{p_3''} \ar[dr]_{f_2} & S_3 \ar@{-->}[dr]^{f_3} & S_1 \ar[d]^{p_3'}\\
      & S_2 \ar[r]_{p_3''} & S_3 \\
      }
  \end{equation*}
\end{proof}

\begin{proposition}
  \label{prop:subsystem_antisymmetry_vec}
  Restrict attention to $f: S \to S$ being a (not necessarily linear) bijection on a vector space $S$, and require that the subsystem projection $p: S \to B$ for each subsystem $(g,p)$ of $f$ is a linear surjection.
  In this case, the relation ``is a subsystem of'' is also antisymmetric up to conjugacy by linear isomorphisms.
\end{proposition}

As a result, data feedback loops are confined to happen within a given subsystem.

\begin{proof}
  Suppose that $(g_2,p_2)$ is a subsystem of $g_1:B_1 \to B_1$, and that $(g_1,p_1)$ is a subsystem of $g_2: B_2\to B_2$, so that we have the commutative diagram
  \begin{equation*}
    \xymatrix{
      B_1 \ar[r]^{g_1} \ar[d]_{p_2} & B_1 \ar[d]^{p_2}\\
      B_2 \ar[r]^{g_2} \ar[d]_{p_1} & B_2 \ar[d]^{p_1} \\
      B_1 \ar[r]^{g_1} & B_1 \\
      }
  \end{equation*}
  Since $p_1$ and $p_2$ are surjective linear maps, this means that $(p_1 \circ p_2): B_1 \to B_1$ is a linear surjection.
  Since it also evidently preserves dimension, it must be a linear isomorphism.
  Because both $p_1$ and $p_2$ are surjective,
  this implies that both must also be injective.
  Hence both $p_1$ and $p_2$ must also be linear isomorphisms,
  which establishes that $g_2 = p_2 \circ g_1 \circ p_2^{-1}$ and $g_1 = p_1 \circ g_2 \circ p_1^{-1}$ as claimed.
\end{proof}

\begin{example}
  There is no function $h$ that will make the diagram below commute
  \begin{equation*}
    \xymatrix{
      \mathbb{Z}^2 \ar[d]_{\id} \ar[r]^f & \mathbb{Z}^2 \ar[d]^{\id} \\
      \mathbb{Z}^2 \ar[r]_{h} & \mathbb{Z}^2 \\
      \mathbb{Z}^2 \ar[u]^{\id} \ar[r]_g & \mathbb{Z}^2 \ar[u]_{\id} \\
      }
  \end{equation*}
  where
  \begin{equation*}
    f(x,y) = (x,1-x),
  \end{equation*}
  and
  \begin{equation*}
    g(x,y)=(y,y).
  \end{equation*}

There is also no function $h$ that will make the diagram below commute
\begin{equation*}
    \xymatrix{
      \mathbb{Z}^2 \ar[d]_{\pr_1} \ar[r]^f & \mathbb{Z} \ar[d]^{\id} \\
      \mathbb{Z} \ar[r]_{h} & \mathbb{Z} \\
      \mathbb{Z}^2 \ar[u]^{\pr_2} \ar[r]_g & \mathbb{Z} \ar[u]_{\id} \\
      }
  \end{equation*}
  where
  \begin{equation*}
    f(x,y) = 1-x,
  \end{equation*}
  and
  \begin{equation*}
    g(x,y)=y.
  \end{equation*}
  
\end{example}

Suppose that $f: S \to S$ is a dynamical system in which $S$ is a vector space and the subsystem projections are all linear surjections, as required by Proposition \ref{prop:subsystem_antisymmetry_vec}.
Let $(\col{B},\le)$ be the collection of all subsystems of $f$, with the partial order established by Lemma \ref{lem:subsystem_preorder} and Proposition \ref{prop:subsystem_antisymmetry_vec}.  Each element of $\col{B}$ is a pair $(g_B,p_B)$ where $g_B: B \to B$ is a bijection and $p_B: S \to B$.  For brevity, if $g_1$ is a subsystem of $g_2$, which is to say that there is a $p_{1,2}: B_2 \to B_1$ such that $p_1 = p_{1,2} \circ p_2$, we write $(g_1,p_1) \le (g_2,p_2)$.

\begin{definition}
\label{df:subsystem_sheaf}
  Define the \emph{sheaf $\shf{F}_f$ of subsystems of $f$} according to the following recipe:
  \begin{description}
  \item[Stalks] $\shf{F}_f((g_B,p_B)) := B$, and
  \item[Restrictions] $\shf{F}_f((g_1,p_1) \le (g_2,p_2)) := p_{1,2}$.
  \end{description}
\end{definition}

Even if the subsystem projections are not linear surjections, the Alexandrov topology on the subsystem preorder bundles together all collections of subsystems that participate in cycles.  Without the conclusion of Proposition \ref{prop:subsystem_antisymmetry_vec}, the stalks of $\shf{F}_f$ are not necessarily well defined, since there is no guarantee that the subsystems of a given cycle have the same state spaces.

\begin{lemma}
  For a dynamical system $f: S \to S$, the space of global sections of $\shf{F}_f$ is precisely $S$.
\end{lemma}
\begin{proof}
  First of all, notice that $\id_S : S \to S$ meets the criteria for a subsystem.  We merely need to verify that the definition of global sections for $\shf{F}_f$ doesn't conflict with this.  The space of assignments for $\shf{F}_f$ is
  \begin{equation*}
    \bigoplus_{p: S\to B \text{ subsystem}} \shf{F}_f(p) = \bigoplus_{p: S\to B \text{ subsystem}} B.
  \end{equation*}
  Suppose that we have a global section $s$.  On the other hand, if $(g_B,p_B) \le (f,\id_S)$, then
  \begin{equation*}
    \left(\shf{F}_f((g_B,p_B) \le (f,\id_S))\right)(s(S)) = p_B(s(S)) = s(B).
  \end{equation*}
  Therefore, the value of $s$ on the subsystem $\id_S:S\to S$ determines the values of $s$ on every other subsystem.
\end{proof}

\begin{proposition}
\label{prop:subsystem_sheaf}
A dynamical system $f: S \to S$ induces an endomorphism on the sheaf of all subsystems, and for which the induced map on global sections is $f$.
\end{proposition}
\begin{proof}
  This follows immediately from the definition, as soon as we notice that for a subsystem $p: S\to B$, the $g$ map guaranteed by the definition is the corresponding component map for the sheaf morphism.
\end{proof}

In short, a multi-scale discrete dynamical system can be encoded as component dynamical systems on some (or all) of the stalks of a sheaf $\shf{S}$ via self maps $f_x:\shf{S}(x) \to \shf{S}(x)$.  One may also consider the action of different semigroups on stalks to model continuous dynamical systems.

We are now ready to establish the main result of this section, which relates the sheaf of subsystems of a DSEM to its graph representation.
As we have seen in Example \ref{eg:logic_oscillator}, feedback loops in the DSEM graph must be confined to being entirely within a subsystem.
Because we can collapse all feedback loops in an arbitrary directed graph to obtain an acyclic graph, we will assume that the DSEM graph is acyclic without loss of generality.

The key insight is that if we select a given variable in the DSEM, any subsystem containing that variable must also contain every variable that can impact its value.
Any variable with a directed path leading to our variable of interest will therefore need to be included in the subsystem.

\begin{definition}
    In a directed graph $G = (V,E)$ an \emph{in-closed} subset $I \subseteq V$ is a set of vertices such that if $v \in I$, then if $e = (w,v) \in E$, then $w \in I$. 
\end{definition}

\begin{lemma}
\label{lem:dsem_subsystem}
    If a dynamical system is defined by a DSEM, every in-closed subset of variables is a subsystem.
\end{lemma}
\begin{proof}
    Suppose that $I$ is a in-closed subset of variables in a DSEM on a directed graph $G$.  
    If $v \in I$ then all of the dependencies of $v$ are also in $I$,
    so the next timestep of $v$ can be predicted from the variables in $I$.
    Therefore, projecting out just the variables in $I$ from the set of all variables will result in a new dynamical update map when restricted to $I$.
\end{proof}

As a consequence of Lemma \ref{lem:dsem_subsystem}, we have the following result that explains why modeling with DSEM is a good idea.

\begin{corollary}
\label{cor:acyclic_dsem_subsystems}
    If a dynamical system is defined by a DSEM on a partially ordered set, then the Alexandrov topology of the dual order is a subspace of the base space topology of its subsystem sheaf.
\end{corollary}

Corollary \ref{cor:acyclic_dsem_subsystems} does not establish that the Alexandrov topology of the dual order of the DSEM \emph{is} the subsystem sheaf.  This is because if the original variables in the DSEM are chosen coarsely, there may be additional subsystems that are ``hidden'' within them.  These hidden subsystems will be present in the subsystem sheaf, but will not correspond to distinct in-closed subsets of the DSEM graph.

\section{Subsystems of the Bering Sea system}
\label{sec:bering_sheaf_subsystems}

\begin{figure}
    \centering
    \includegraphics[width=\linewidth]{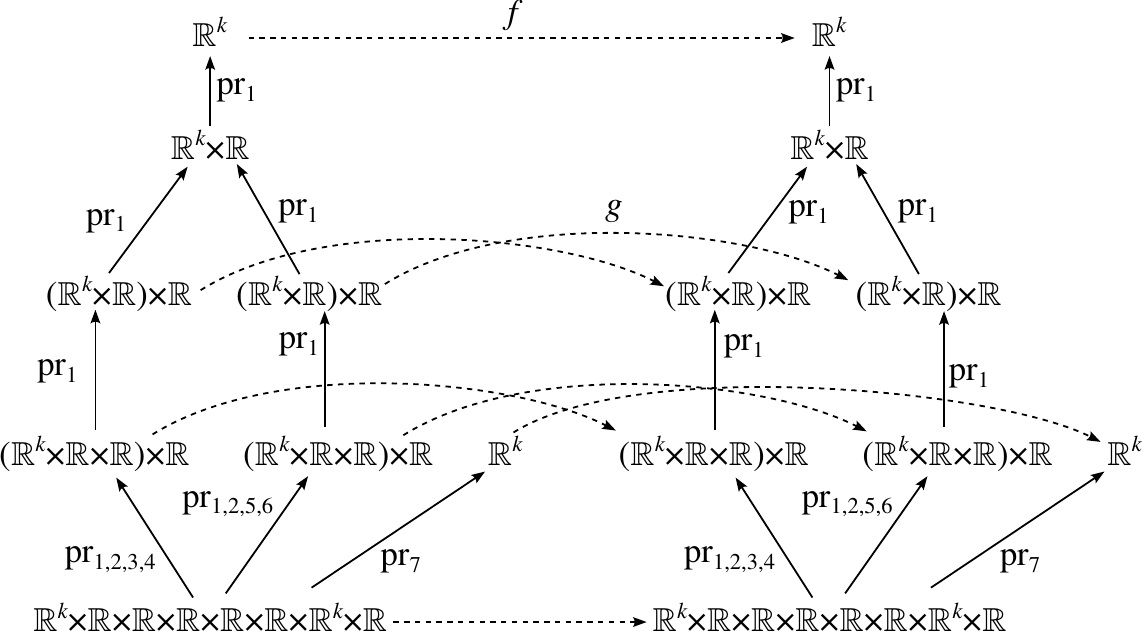}
    \caption{Sheaf of subsystems for the Bering Sea example.  Solid arrows are the subsystem projection maps; dashed arrows are the dynamical system state update maps.  Maps $f$ and $g$ are explained in the text.}
    \label{fig:bering_sea_subsystems}
\end{figure}

Figure \ref{fig:bering_sea_subsystems} shows the sheaf of subsystems for the Bering Sea example, with the stalks organized in the same way as shown in Figure \ref{fig:bering_sea_sheaf}.  

The function $f$ performs an $\AR(k)$ update:
\begin{equation*}
    f(x_1, \dotsc, x_k) =\left( x_2, \dotsc, x_k, \sum_{i=0}^{k-1} a_i x_{k-i}\right),
\end{equation*}
while the function $g$ performs the dynamical update for the subsystem containing the $\mathit{Krill}$ variables:
\begin{equation*}
    g(x_1, \dotsc, x_k, y, z) = \left(x_2, \dotsc, x_k, \sum_{i=0}^{k-1} a_i x_{k-i}, y + c x_k, z + d y\right). 
\end{equation*}
Notice how $f$ is obtained from $g$ by projecting out the first $k$ components, in accordance with the commutativity of Figure \ref{fig:bering_sea_subsystems}.

Although Figures \ref{fig:schematic_sheafify}(d) (with modifications to support autoregressive timeseries), \ref{fig:bering_sea_sheaf}, and \ref{fig:bering_sea_subsystems} represent different sheaves, they all represent the same dynamical system.  
Consequently, the global sections of these three sheaves are different but are in a natural bijective correspondence.
The three sheaves offer three distinct perspectives, with increasing granularity,
\begin{description}
    \item[Definition \ref{df:subsystem_sheaf}: Figure \ref{fig:bering_sea_subsystems}] Stalks are nested collections of dynamically related variables, each represented by sliding windows of timeseries,
    \item[Definition \ref{df:netlist_sheaf}: Figure \ref{fig:schematic_sheafify}(d)] Each variable is an entire timeseries and appears alone in at least one stalk, and
    \item[Definition \ref{df:netlist_sheaf_individual}: Figure \ref{fig:bering_sea_sheaf}] Each observation (a timestep for a single variable) appears alone in at least one stalk.
\end{description}
With this perspective, the boundaries between subsystems are easily seen in Figure \ref{fig:bering_sea_sheaf}: those restriction maps that are identity maps from parts to nets are those that cross subsystem boundaries.  The variables at the heads of any identity maps in Figure \ref{fig:bering_sea_sheaf} are those that are removed by the subsystem projections involved.  Moreover, the state spaces arise as one time step of the space of local sections over each subsystem, once cut.

\section{Conclusion}
\label{sec:conclusion}

In this chapter, we have demonstrated how the general framework of sheaf modeling applies to several composite dynamical systems, including an ecological model of the Bering Sea and a dynamical model of low-level computer software.
Sheaf modeling provides a coherent mathematical framework for studying the complicated interaction of various dynamical subsystems that together determine a larger system.
The guiding principles of sheaf modeling are that
\begin{itemize}
    \item a sheaf represents a hypothesis about how variables will interact,
    \item a non-global assignment represents the observations collected on the variables in its support,
    \item minimizing consistency radius predicts values of the variables that were not observed, and
    \item the minimal consistency radius is a measure of the consistency between the observations and the hypothesis.
\end{itemize}
This chapter shows that when a dynamical system is described by a DSEM, there are three sheaves that provide increasingly granular data about the interactions between variables:
\begin{enumerate}
    \item the sheaf of subsystems (Definition \ref{df:subsystem_sheaf}),
    \item the netlist sheaf with timeseries as stalks (Definition \ref{df:netlist_sheaf}), and
    \item the netlist sheaf with additional stalks for individual observations (Definition \ref{df:netlist_sheaf_individual}).    
\end{enumerate}
With these three sheaves in hand, a system modeler can apply the guiding principles above to measure how well their model fits observational data.
The sheaf encodings allow the modeler to perform a variety of standard inferences (e.g. forward prediction, backward prediction, regression, and missing-data imputation) using a unified framework. 
The sheaf modeling framework easily supports hybrid versions, for instance performing simultaneous forward and backward predictions, or simultaneously performing regression and prediction.
Since the sheaf framework measures the fit between observations and the model, the modeler can assess their confidence in these inference tasks.

It remains future work to compare estimates of uncertainty computed by the DSEM (appearing in the $\mathbf{V}$ and $\mathbf{E}$ matrices) to the consistency radius of the corresponding sheaf.  In particular, it seems possible to view consistency radius as a test statistic for the distributional model posited by the DSEM.  Indeed, Equation \eqref{eq:consistency_radius} is strikingly close to the log likelihood if the distributions of measurement errors are assumed to follow an exponential model.  If this is true, then it should be possible to lift the sheaf modeling discipline described here into a standard statistical hypothesis testing framework.  

\section*{Acknowledgments}

The linear regression example in Section \ref{sec:netlist_sheaf} is due to Donna Dietz.

This article is based upon work supported by the Office of Naval Research (ONR) under Contract Nos. N00014-15-1-2090 and N00014-18-1-2541, the Defense Advanced Research Projects Agency (DARPA) SafeDocs program under contract HR001119C0072, and the MITRE Corporation's Independent Research and Development (IR\&{}D) Program.  Any opinions, findings and conclusions or recommendations expressed in this article are those of the authors and do not necessarily reflect the views of ONR, DARPA, or MITRE.

\bibliographystyle{plainnat}
\bibliography{netlist_sheaf_bib,netlist-sheaf-thorson}

\end{document}